\documentclass{amsart}
\usepackage[utf8]{inputenc}
\usepackage[margin=1in,left=1in]{geometry}
 \usepackage{mathtools}
\usepackage{setspace}
\usepackage{graphicx}
\usepackage{amssymb, amsmath, amsthm, graphics}
\usepackage{mathabx,epsfig}
\usepackage[mathscr]{euscript}
\usepackage{tikz}
\usepackage{pigpen}
\usepackage{bm}
\usepackage{mathtools}
\usetikzlibrary{calc}
\usetikzlibrary{matrix}
\usepackage{latexsym}
\usepackage[all]{xy}
\usepackage{color}

\usepackage{bbold}

\usepackage{tikz}
\input xy
\xyoption{all}
\usepackage{euscript}
\usepackage{multirow}
\usepackage{etex, pictexwd,dcpic}
\usetikzlibrary{positioning}
\usetikzlibrary{shapes.geometric}
\usetikzlibrary{shapes.misc}
\usetikzlibrary{calc}
\usetikzlibrary{positioning}

 \usepackage{enumerate}

\usepackage{pgflibraryarrows}
\usepackage{pgflibrarysnakes}

\usetikzlibrary{trees} 
\usetikzlibrary[trees] 

\newtheorem{theorem}{Theorem}
\newtheorem{definition}[theorem]{Definition}
\newtheorem{lemma}[theorem]{Lemma}

\newtheorem{proposition}[theorem]{Proposition}

\newtheorem{remark}{Remark}

\newtheorem{example}{Example}

\numberwithin{equation}{section}

\renewcommand{\(}{\begin{equation*}}
\renewcommand{\)}{\end{equation*}}
\newcommand{\bea}{\begin{eqnarray*}}
\newcommand{\eea}{\end{eqnarray*}}

\renewcommand{\a}{\alpha}
\renewcommand{\b}{\beta}

\def\endofproof {\hfill{$\Box$}\\}

\def\H{\ensuremath{\ES{H}}}

\def\H{\ensuremath{\ES{H}}}

\newcommand{\beq}{\begin{equation}}
\newcommand{\eeq}{\end{equation}}

\newcommand{\into}{\hookrightarrow}

\newcommand{\op}[1]{\ensuremath{\operatorname{#1}}}

\newcommand{\ES}[1]{\ensuremath{\EuScript{#1}}}

\newcommand{\theproof}{\noindent {\bf Proof.\ }}

\numberwithin{equation}{section}

\renewcommand{\(}{\begin{equation}}
\renewcommand{\)}{\end{equation}}

\def\ch{{\rm  ch}}

\def\1{{\bf 1}}

\def\<{\langle}
\def\>{\rangle}

\def\a{\alpha}
\def\b{\beta}

\numberwithin{equation}{section}


\newcommand{\R}{\ensuremath{\mathbb R}}
\newcommand{\RR}{\ensuremath{\mathbb R}}

\newcommand{\ZZ}{\ensuremath{\mathbb Z}}
\newcommand{\Z}{\ensuremath{\mathbb Z}}

\newcommand{\BB}{\ensuremath{\mathbf B}}

\newcommand{\sh}{\ensuremath{\mathscr{S}\mathrm{h}}}

\newcommand{\map}{\mathrm{Map}}

\begin{document}

\title{Higher-twisted periodic smooth Deligne cohomology}

 \author{Daniel Grady and Hisham Sati
 \\
  }

\maketitle

\begin{abstract} 
Degree one twisting of  Deligne cohomology, as a differential refinement of 
integral cohomology, was established in previous work. 
Here we consider higher degree twists.  The Rham complex, hence de 
Rham cohomology,  admits twists of any odd degree. However, in order to 
 consider twists of integral cohomology we need a periodic version.
Combining the periodic versions of both ingredients leads us to introduce 
a periodic form of Deligne cohomology. We demonstrate that this theory indeed admits 
a twist by a gerbe of any odd degree. We present the main properties of 
the new theory and illustrate its use with examples and computations, mainly 
via a corresponding twisted differential Atiyah-Hirzebruch spectral sequence. 
  \end{abstract}

 \tableofcontents
 
\section{Introduction}

There has been a lot of recent activity on modifying generalized cohomology 
theories to include twists and geometric refinements, in order to account 
for automorphisms  and  include geometric data. 
Twisted differential generalized cohomology theories are established at the general axiomatic 
level  \cite{BN}. However, working out these theories 
explicitly is in practice not a straightforward task. Twisting the simplest case of 
a differential cohomology theory, namely Deligne cohomology, proved to be
nontrivial \cite{GS4} and is closely related to interesting constructions in 
algebraic geometry, namely taking coefficients in variations of mixed Hodge structures 
(see \cite{CH}\cite{H2}). Even at the topological level, while twisting of generalized 
cohomology theories is axiomatically well-established \cite{MS}\cite{Units}\cite{ABG}, 
spelling out explicit constructions requires considerable work (see 
\cite{ABG}\cite{SW}\cite{LSW} for recent illustrations). The goal of this paper 
is to generalize the degree one twists of Deligne cohomology from \cite{GS4}  to include twists 
of higher degrees. These will be in the form of higher gerbes, or $n$-bundles, 
with connections (see \cite{Cech}\cite{SSS3}\cite{Urs}\cite{FSS1}\cite{FSS2} 
for constructions and related applications). 

\medskip
Deligne cohomology (see \cite{Del}\cite{Be}\cite{Gi}\cite{Ja}\cite{EV}\cite{Ga1}) is 
a differential refinement of ordinary, i.e. integral, cohomology. As such it has various 
realizations (see 
 \cite{Del}\cite{CS}\cite{Ga1}\cite{Bry}\cite{DL}\cite{HS}\cite{BKS}
 \cite{BB}\cite{Urs}), 
 which  are (expected to be) equivalent   (see \cite{SSu}\cite{BS}). 
Consider the sheaf of chain complexes associated with 
real-valued differential forms
\footnote{This is sometimes also denoted $\mathbb{Z}_{\mathcal D}^\infty(n)$ or $\mathbb{Z}(n)_{\mathcal D}^\infty$.
We are in the smooth setting throughout, so we will not need extra decorations.}
\(
\label{Delcplx}
{\mathcal D}(n):=\big(
\xymatrix@C=.8cm{ 
\hdots \ar[r] & 0\ar[r] & \underline{\ZZ}\ar[r]  & \Omega^0\ar[r]^{d} &
 \Omega^1\ar[r]^{d} & \hdots \ar[r]^{d} & \Omega^{n-1}
}\big)\;,
\)
where we place differential $(n-1)$-forms in degree $0$ and locally constant integer-valued 
functions in degree $n$. Given a smooth manifold $M$, the Deligne cohomology group of 
degree $n$ is defined to be the sheaf  (hyper-)cohomology group 
\footnote{This would be $H^n(M; {\mathcal D}(n))$ if we use the opposite convention.
However, the one we use is positively graded, hence better adapted for stacks.}
$\widehat{H}^n(M;\ZZ):=H^0(M;{\mathcal D}(n))$. 
 {\v C}ech resolutions allow for explicit calculation of these groups.  
 If $\{U_{\alpha}\}$ is a good open cover of $M$, then one can form
  the {\v C}ech-Deligne double complex
(see \cite{BT}\cite{Bry})
\begin{equation}
\hspace{-1mm}
{\small
\label{double complex} 
\xymatrix@C=1em{
 \prod \limits_{\a_0, \cdots, \a_n}\underline{\ZZ}(U_{\alpha_0\hdots \alpha_{n}})
 \ar[r] & 
 \prod \limits_{\a_0, \cdots, \a_n} \Omega^0(U_{\alpha_0\hdots \alpha_{n}})\ar[r]^-{d} & 
  \prod \limits_{\a_0, \cdots, \a_n}\Omega^1(U_{\alpha_0\hdots \alpha_{n}})\ar[r]^-{d} &
  \hdots \ar[r]^-{d} &
   \prod \limits_{\a_0, \cdots, \a_n} \Omega^{n-1}(U_{\alpha_0\hdots \alpha_{n}})
\\
 \prod \limits_{\a_0, \cdots, \a_{n-1}}\underline{\ZZ}(U_{\alpha_0\hdots \alpha_{n-1}})
 \ar[r] \ar[u]^-{ (-1)^{n-1}\delta} &  \prod \limits_{\a_0, \cdots, \a_{n-1}} 
 \Omega^0(U_{\alpha_0\hdots \alpha_{n-1}})
\ar[r]^-{d} \ar[u]^-{(-1)^{n-1}\delta}& 
 \prod \limits_{\a_0, \cdots, \a_{n-1}}\Omega^1(U_{\alpha_0\hdots \alpha_{n-1}})
 \ar[r]^-{d}\ar[u]^-{(-1)^{n-1}\delta} & 
\hdots \ar[r]^-{d} &  \prod \limits_{\a_0, \cdots, \a_{n-1}} \Omega^{n-1}(
U_{\alpha_0\hdots \alpha_{n-1}})\ar[u]^-{(-1)^{n-1}\delta}
\\
\vdots \ar[u]^-{(-1)^{n-2}\delta}& \vdots \ar[u]^-{(-1)^{n-2}\delta}& \vdots \ar[u]^-{(-1)^{n-2}\delta}
&  & \vdots \ar[u]^-{(-1)^{n-2}\delta}
\\
  \prod \limits_{\a_0, \a_1}\underline{\ZZ}(U_{\alpha_0\alpha_1})\ar[r]\ar[u]^-{-\delta} & 
\prod \limits_{\a_0, \a_1} \Omega^0(U_{\alpha_0\alpha_1})\ar[r]^-{d}\ar[u]^-{-\delta} &
\prod \limits_{\a_0, \a_1}  \Omega^1(U_{\alpha_0\alpha_1})\ar[r]^-{d}\ar[u]^-{-\delta} & 
  \hdots \ar[r]^-{d} &\prod \limits_{\a_0, \a_1} \Omega^{n-1}(U_{\alpha_0\alpha_1})\ar[u]^-{-\delta}
\\
\prod \limits_{\a_0}\underline{\ZZ}(U_{\alpha_0})\ar[r] \ar[u]^-{\delta} & 
\prod \limits_{\a_0}\Omega^0(U_{\alpha_0})\ar[r]^-{d}\ar[u]^-{\delta} & 
\prod \limits_{\a_0}\Omega^1(U_{\alpha_0})\ar[r]^-{d}\ar[u]^-{\delta} & \hdots \ar[r]^-{d} &
\prod \limits_{\a_0} \Omega^{n-1}(U_{\alpha_0})\ar[u]^-{\delta}\;,
}
}
\end{equation}
where $U_{\alpha_0\alpha_1\hdots \alpha_k}$ denotes the $k$-fold intersection
$U_{\a_0} \cap U_{\a_1} \cap \cdots \cap U_{\a_k}$. 
With $d$ and $\delta$ 
 the de Rham and {\v C}ech differentials, respectively, acting on elements of degree $p$,
the total operator on the double complex is the {\v C}ech-Deligne operator 
$D:=d+(-1)^p\delta$.
The sheaf cohomology group $H^0(M;{\mathcal D}(n))$ can then be identified with 
the group of diagonal 
elements $\eta_{{}_{k,k}}$ in the double complex which are $D$-closed,
$D\eta_{{}_{k,k}}=0$, modulo those which are $D$-exact.

\medskip
At a more general level, and from a homotopy theory point of view, given a spectrum 
$E$ the canonical data for the corresponding differential theory is comprised of the 
following (see \cite[Example 4.49]{Bu}).  Let $\H$ be the Eilenberg-MacLane functor and take 
$A:=\pi_* E \otimes \RR$ the `realified' coefficients of the theory.
 Let $c: E \to \H(A)$ be the map uniquely determined up to homotopy 
 such that it induces the map `realifying' the coefficients 
 $\pi_*(E) \to \pi_*(A) \cong \pi_*(E) \otimes \RR$, $x \mapsto x \otimes 1$. 
 Indeed, for $E=\H(\Z)$ the integral Eilenberg-MacLane spectrum, 
 $A=\Z \otimes \RR \cong \R$, and $c: \H(\Z) \to \H(\RR)$ 
uniquely determined by $\Z \to \RR\cong \Z \otimes \RR$, $x \mapsto x \otimes 1$. 
This data determines a differential extension of $\H \ZZ$, which in 
degree $n$ takes the form $(\Sigma^n\H\Z, \RR, c)$.
Applying the Eilenberg-MacLane functor $\H$ to the Deligne complex ${\mathcal D}(n)$ from 
expression \eqref{Delcplx} gives a natural equivalence of  differential  
spectra (see \cite{Bu})
$$
\H({\mathcal D}(n)) \cong (\Sigma^n\H\Z, \RR, c)\;.
$$
%

\medskip
What would it mean to twist Deligne cohomology? At a practical level, Deligne cohomology is 
a fusion of de Rham cohomology and {\v C}ech cohomology. One could naively 
consider twists of de Rham cohomology, made periodic, as well as twists of 
{\v C{ech cohomology,
separately. 
Now, the {\v C}ech-Deligne operator  $D=d + (-1)^p\delta$ acting on the {\v C}ech-Deligne 
double complex is a combination of the de Rham differential $d$, acting in the direction of the de Rham
complex, and the {\v C}ech differential $\delta$, acting in the direction of the {\v C}ech 
complex (see again \cite{BT}\cite{Bry} for extensive 
discussions). One could then consider two situations (or stages) in attempting to 
twist the  operator $D$: First, adding a closed differential form $H$, i.e. modifying 
only the form part, leading to twisting of the de Rham differential $d \leadsto d_H$. 
Second, adding a differential cohomology class $\hat{h}$, i.e. 
modifying both parts of the whole {\v C}ech-de Rham differential, i.e. $d \leadsto d_H$ and 
$\delta \leadsto \delta_{B}$ where $(H, B)$ are appropriate (de Rham, {\v C}ech)-components 
of $\hat{h}$. We will see that it will not quite work that way; nevertheless, this turns out to be a good heuristic 
to keep in mind.  In contrast to the de Rham differential, twisting the {\v C}ech differential would be 
quite complicated due to  higher and higher local transition data. 
Indeed, {\v C}ech-de Rham cohomology, twisted by a degree three form (without altering 
the {\v C}ech direction) is used in \cite{GT2} for describing a finite-dimensional model 
of twisted K-theory. Note that twisted de Rham cohomology of a space can be described 
via the untwisted  cohomology of a corresponding stack \cite{BSS}.

\medskip
Which degrees should the twists $H$ or $\hat{h}$ have? 
The twists of the Deligne complex, a priori naturally arise in degree one \cite{GS4}. 
Note that for the underlying topological theory, a representation of the fundamental 
group $\pi_1(X)$ of a space $X$ on $\op{Aut}(\Z)\cong \Z/2$
gives $\ZZ$ the structure of a module over the group ring $\Z[\pi_1(X)]$,
which is used in \cite{BFGM} to describe $\pi_1(X)$-twisted integral cohomology.
On the other hand, one can twist the de Rham complex by differential forms of any odd degree, 
not just degree one  (see \cite{RW}\cite{BCMMS}\cite{Tel}\cite{higher}\cite{tcu}\cite{MW}).  
 At first glance, this might appear to give an inherent incompatibility of twisted de Rham cohomology 
 and twisted integral cohomology. However, if one takes a closer look,
 one  realizes that twisted de Rham cohomology is really about $\ZZ/2$-graded, periodic de Rham 
  cohomology. Thus, one does not expect compatibility with integral cohomology, but rather with 
  \emph{periodic} integral cohomology. Hence we consider twists of the latter theory in 
  Section \ref{Sec Per Z}.   This then paves part of the way for us to go towards a general twisted 
  Deligne cohomology. However, as both ingredients, namely de Rham and integral cohomology, were 
  made periodic, we define a periodic version of Deligne cohomology in Section \ref{Sec Per Del}. 
  We characterize its main properties via sheaf cohomology and differential spectra, including 
  the ring structure arising from the Deligne-Beilinson cup product 
\cite{Del}\cite{Be} (see \cite{FSS1}\cite{FSS2}). 
Periodic integral differential cohomology groups 
$\widehat{H}(X; \ZZ[u, u^{-1}])$ have been considered 
from an index theoretic point of view briefly  in \cite{Lo}\cite[Sec. 8.4]{FL}.

\medskip
Having defined the appropriate starting point for the twisting of Deligne cohomology, namely the 
periodic version, we discuss the twists of periodic Deligne cohomology in Section \ref{Sec Higher}. 
We approach twisting of periodic  Deligne cohomology using simplical presheaves and 
smooth stacks \cite{Cech}\cite{FSS1}\cite{HQ}\cite{FSS2}\cite{Urs}, as we did 
in \cite{GS4}. This approach is very well-suited to the higher twists and 
 allows for the use of powerful algebraic machinery. We will show that the twists indeed 
 refine the twists of both integral cohomology and  the de Rham complex. 
Smooth stacks will arise naturally in twisting periodic Deligne cohomology. Just as we can 
 twist periodic integral cohomology by odd degree singular cocycles (Section \ref{Sec Per Z}), 
 we will see that periodic Deligne cohomology can be twisted by higher gerbes of odd degree
 (Section \ref{Sec Higher}). The appearance of gerbes naturally leads us into the world of smooth 
 stacks, and we will find it useful to recall some of the constructions in this setting 
 (see \cite{Bry}\cite{Cech}\cite{FSS1}\cite{FSS2}).
This requires us to understand in detail exactly what we mean by twisting a periodic 
differential cohomology theory. We  give a characterization of the twists 
via moduli stacks of higher bundles with connections, but keeping technical 
matters to a minimum.

\medskip
As discussed above, we can twist Deligne cohomology by gerbes of odd degree. It is interesting 
to see where the gerbe data appears in defining the twisted theory. In fact, as observed in 
\cite{BN}, a crucial ingredient in defining twisted differential theories is the analogue 
of the de Rham isomorphism theorem for twisted cohomology. In the untwisted case, 
recall that the locally constant sheaf $\underline{\RR}$ admits an acyclic resolution via 
the de Rham complex
\(\label{de Rham resolution}
\xymatrix{
\underline{\RR}\ \ar@{^{(}->}[r] & \Omega^0\ar[r]^{d} & \Omega^1\ar[r]^{d} & 
\Omega^2 \ar[r] & \hdots \;,
}
\)
and the de Rham Theorem is manifestly a corollary of this fact. Indeed, for a smooth 
manifold $M$, the sheaf cohomology $H^\ast(M,\underline{\RR})$ can be calculated 
both as {\v C}ech cohomology and via this resolution. The isomorphism between 
singular and {\v C}ech cohomology then recovers
 de Rham's classical theorem. 


\medskip
Just as multiplicative cohomology theories have topological spaces of twists (the Picard spaces)
\cite{MQRT}
\cite{MS}\cite{Units}, differential refinements of such theories have \emph{smooth stacks} 
of twists. Indeed,  the stack of twists $\widehat{{\rm Tw}}_{\widehat{\mathscr{R}}}$ 
for any differentially refined cohomology theory $\widehat{\mathscr{R}}=(\mathscr{R},c,A)$ 
was introduced in \cite{BN}. This was defined by the pullback
\footnote{
We will be dealing with $(\infty, 1)$-categories, so that whenever we talk about 
pullbacks, pushouts, or any other universal construction, we mean it in the 
$(\infty, 1)$-sense, i.e., up to higher coherence  homotopy. Whenever we draw 
a diagram, it should be understood, whether we draw it explicitly or not, that there are
homotopies and higher homotopies involved, up to the appropriate degree!}
   (in the notation of \cite{GS5})
\(
\label{Tw-diag}
\xymatrix@R=1.5em{
\widehat{\rm Tw}_{\widehat{\mathscr{R}}}\ar[rr]\ar[d] && 
{\rm Pic}_{\mathscr{R}}^{\rm form}\ar[d]
\\
{\rm Pic}_{\mathscr{R}}^{\rm top}\ar[rr] &&
 {\rm Pic}_{\mathscr{R}}^{\rm dR}\;,
}
\)
where 
\begin{itemize}

\vspace{-2mm}
\item ${\rm Pic}_{\mathscr{R}}^{\rm top}$ is the ordinary Picard $\infty$-groupoid 
of twists for the ring spectrum $\mathscr{R}$, embedded as a constant smooth stack, 

\item ${\rm Pic}_{\mathscr{R}}^{\rm dR}$ 
is the Picard stack of \emph{sheaves} of invertible module spectra over the smash product 
$\mathscr{R}\wedge \H\RR$ (embedded as a constant sheaf of spectra), and 

\item ${\rm Pic}_{\mathscr{R}}^{\rm form}$ is the smooth stack which (after evaluation 
on a smooth manifold $M$) comes as the nerve of the groupoid whose objects are weakly 
locally constant, K-flat, invertible modules over $\Omega^*(-;A)\vert_{M}$ 
(see \cite{BN} for details). 
\end{itemize}

\noindent An element of the pullback \eqref{Tw-diag} can be identified with a triple 
$\widehat{\mathscr{R}}_{\hat{\tau}}=(\mathscr{R}_{\tau},t,{\mathcal L})$, 
where $\mathscr{R}_{\tau}$ is an underlying twisted cohomology theory 
with a topological twist $\tau$, ${\mathcal L}$ is an invertible module over 
$\Omega^*(-;A)$ and $t$ is an equivalence
$$
t:\mathscr{R}_{\tau}\wedge \H\RR\overset{\simeq}{\longrightarrow} \H({\mathcal L})\;,
$$
exhibiting a twisted de Rham theorem. This stack will be important in 
identifying the twists for periodic Deligne cohomology in Section \ref{Sec Higher}.

\medskip
The situation is summarized in the 
following tables. The first one captures the untwisted case

\medskip
\begin{center}
\begin{tabular}{|c||c|c|}
\hline
 {\bf Untwisted cohomology} & {\it Ordinary} & {\it Periodic} 
\\
\hline\hline
{\it Underlying theory} & Locally constant sheaf $\underline{\RR}$ 
& Sheaf of graded algebras $\underline{\RR}[u, u^{-1}]$ 
\\
\hline
{\it de Rham complex} & Ordinary de Rham complex $\Omega^*$ 
& Periodic complex $\Omega^*[u, u^{-1}]$ 
\\
\hline
\end{tabular}
\end{center}

\medskip
In \cite{GS5} we highlighted the close analogies between twisted spectra 
and line bundles, in that twisted differential spectra are closely related to 
bundles of spectra equipped with a flat connection. Here, we wish to replace 
the  notions in the above table with the twisted analogues as follows.

\medskip
\begin{center}
\hspace{-2mm}
\begin{tabular}{|c||c|c|}
\hline
 {\bf Twisted cohomology} & {\it Ordinary} & {\it Periodic} 
 \\
 \hline
 \hline
 {\it Underlying theory}  & Locally constant sheaf ${\mathcal L}$ 
 & Sheaf of DGA-modules  $ {\mathcal L}_{\bullet}$
\\
\hline
{\it Twist degree} & One & Any odd degree
\\
\hline
{\it Geometric twisting object} & Line bundle 
with flat connection $d+H_1$ 
& Gerbe with curvature $H_{2k+1}$
\\
\hline
{\it de Rham complex} &  ($\Omega^* \otimes {\mathcal L}$,  $d+H_1\wedge$)
& ($\Omega^* \otimes {\mathcal L}_{\bullet}$, $d+H_{2k+1}\wedge$)
\\
\hline
\end{tabular}
\end{center}

\medskip
The theories that we consider are related schematically as follows
\(
\label{theories}
\xymatrix{
\fbox{Periodic Deligne $\widehat{H}^*(M; {\mathcal D}[u, u^{-1}])$} 
\ar@{~>}[d]_-{\tiny \begin{array}{cc} {\tiny \rm forget}\\{\rm connection}\end{array}}
\ar@{~>}[rrr]^-{\rm forget}_-{\rm periodicity}
&&&
\fbox{Deligne $\widehat{H}^*(M; \ZZ)=H^0(M; {\mathcal D}(*))$} 
\ar@{~>}[d]^-{\tiny \begin{array}{cc} {\tiny \rm forget}\\{\rm connection}\end{array}} 
\\
\fbox{Periodic integral  ${H}^*(M; \ZZ[u, u^{-1}])$} 
\ar@{~>}[rrr]^-{\rm forget}_-{\rm periodicity}
&&&
\fbox{Integral   $H^*(M; \ZZ)$}
}
\)
where the top row, bottom row, left column and right column represents
geometric theories, topological theories, periodic theories, and 
non-periodic theories, respectively. 
The relations between the corresponding spaces of twists are in turn summarized in 
the schematic diagram 
\(
\label{twist-spaces}
{
\xymatrix{
\BB(\ZZ/2)_\nabla \times \prod_{k >0} \BB^{2k} U(1)_\nabla \; 
\ar@{~>}[d]_{|\;\cdot\;|}  \ar@{~>}[rrr]^-{u=0}
&&& \;
\BB(\ZZ/2)_\nabla \ar@{~>}[d]^-{|\;\cdot\;|} 
\\
K(\ZZ/2, 1) \times \prod_{k>0} K(\ZZ, 2k+1) \ar@{~>}[rrr]^-{u=0}
&&&
K(\ZZ/2, 1)\;.
}
}
\)
Here $|\cdot|$ is geometric realization, which reduces a geometric theory 
down to the corresponding topological theory,  and $\BB(\ZZ/2)_\nabla$ is the 
stack of twists  for Deligne cohomology \cite{GS4}. 
To twist the  theories displayed in the first schematic diagram
\eqref{theories} one would consider maps from the manifold 
$M$ to the corresponding space of twists in 
the second schematic diagram \eqref{twist-spaces}.

\medskip
Explicit {\v C}ech cocycles for Deligne cohomology are described in 
\cite{Ga1}\cite{BM1}\cite{BrM}\cite{Gomi1}.  While we do not do this in full generality
in the twisted case, we do explain how the {\v C}ech cocycle data appear as part
of the trivializing data for twisted periodic Deligne cohomology in 
Remark \ref{rem-CS}.
 This involves Chern-Simons type trivialization 
of  {\v C}ech-Deligne cocycles, packaged succinctly as in \cite{GS5}.
Extensive discussions of such  trivializations relating to Chern-Simons
theory are given in  \cite{BM1}\cite{BrM}\cite{Gomi2}\cite{Fr2}\cite{CJMSW}\cite{Cech}\cite{Wa}
\cite{GPT}\cite{M5}\cite{tert}\cite{FSS2}\cite{E8}\cite{Th}. 
In contrast, cocycles arising from chain complexes
would involve an abstract higher local system resulting from the twists of periodic 
integral cohomology. While this is doable, it does not make the description 
any more transparent in comparison to the description via spectra; hence 
we do not consider it in this paper. 

\medskip
Making use of the general 
constructions in \cite{GS3}\cite{GS5}, 
we then consider the Atiyah-Hirzebruch spectral sequence for 
 twisted periodic integral cohomology as well as for twisted 
periodic Deligne cohomology in Section \ref{Sec-tAHSS}. We  
provide explicit constructions and characterizations in Section 
\ref{tAHSS} and then illustrate the computations  via 
 examples in Section \ref{Sec-Ex}.

\medskip
We note that there are other approaches to studying the {\v C}ech-de Rham double 
complex. Explicit description of cocycles via the cohomology of the total 
operator $D$ of the double complex is provided in \cite{Pi}. 
Using the notion of Cheeger-Simons cochain sparks \cite{CS}, a homological 
machine for the study of secondary geometric invariants called spark complexes is 
described in \cite{HLZ} \cite{HL}. This seems to be an appropriate setting for twisting differential 
cohomology in its incarnation as differential characters. While we do not address this, we expect 
that the resulting twisted versions would be equivalent;  
our approach places the complication in 
the coefficients of the (hyper)cohomology while that more homologically flavored approach would 
place it in the cycles, e.g., via spark complexes.

\section{Twisted periodic integral cohomology}
\label{Sec Per Z}

In this section, we describe the twisted periodic cohomology with both real and integral 
coefficients. This generalizes twists of integral cohomology \cite{MQRT}, also described in 
modern categorical terms in \cite{ABG} and geometrically in \cite{Fr}. This will be a precursor for the 
de Rham theorem needed to define twisted Deligne cohomology. Throughout the remainder 
of this paper, we will follow the $\infty$-categorical treatment of twisted cohomology theories
 \cite{ABG}\cite{Units}\cite{SW}  and their differential refinements 
\cite{BN}\cite{GS5}.

\subsection{Twists via bundles of spectra}
\label{Sec-bunSp} 

The starting point for periodic integral cohomology is the differential graded algebra 
(DGA) $\ZZ[u, u^{-1}]$, 
equipped with the trivial differential. There is a functor
\(
\label{EM-fun}
\H:\mathscr{C}{\rm h}\longrightarrow \mathscr{S}{\rm p}\;,
\)
from the category of unbounded chain complexes to the category of spectra,
 called the Eilenberg-MacLane functor. This functor was defined in \cite{Shi}, where it was 
 shown to exhibit an equivalence between 
$\H\ZZ$-module spectra and differentially graded $\ZZ$-algebras. Applying $\H$ to $\ZZ[u, u^{-1}]$ 
we get a spectrum $\H\ZZ[u,u^{-1}]$ which represents periodic integral cohomology,
 in the sense that 
\(
\label{iso zper}
H^*(X;\ZZ[u, u^{-1}])\cong H^*(X;\ZZ)[u, u^{-1}]\;,
\)
where the right hand side is the graded algebra whose elements are formal Laurent polynomials with 
coefficients in $H^*(X;\ZZ)$ graded by homogeneous degree. The ring structure on the right is induced 
from the cup product structure on $H^*(X;\ZZ)$, while on the left it is induced from the algebra 
structure on $\ZZ[u,u^{-1}]$. This theory is naturally $\ZZ/2$-graded, as we have canonical isomorphisms 
$H^*(X;\ZZ[u,u^{-1}])\cong H^{*+2}(X; \ZZ[u, u^{-1}])$. 
For this reason, we will usually refer to the degree of a class as either \emph{even} or \emph{odd}.

\begin{remark} 
[Action of units on periodic integral cohomology]
Being an $\H\ZZ$-module spectrum, the spectrum $\H\ZZ[u,u^{-1}]$ receives an action by $\H\ZZ$. 
This action manifests itself simply by the action of the cup product in integral cohomology. More precisely, 
given a cohomology class $h\in H^{\ast}(X;\ZZ)$, we can act on a homogeneous Laurent polynomial 
of degree $k$,  $a=\cdots +a_{k+2}u^{-1}+a_k+a_{k-2}u^{1}+\cdots $, via
$$
(h,a)\longmapsto \; \cdots +(ha_{k+2})u^{-1}+(ha_k)+(ha_{k-2})u^{1}+\cdots\;.
$$
Notice, however, that if $h$ has odd degree then the action
would not arise from the units of the spectrum $\H \ZZ[u, u^{-1}]$,
since the powers of $u$ are weighted by even degrees. On the other hand, for an even degree cohomology class
$h$, the action by $h u^{-{\rm deg}(h)}$ preserves the homogenous degree and gives rise to twist 
of $H^*(X;\ZZ[u,u^{-1}])$ on $X$. 
\end{remark}

We now characterize the space of twists of periodic integral cohomology,
the first summand of which is the space of twists in the non-periodic case 
described in \cite{GS4}. 

\begin{proposition}
[Space of twists for periodic integral cohomology]
\label{prop twhz}
The space of twists for periodic integral cohomology is

\vspace{-3mm}
$$
B{\rm GL}_1(\H\ZZ[u,u^{-1}]) \simeq K(\ZZ/2,1)\times \prod_{k>0}K(\ZZ,2k+1)\;.
$$
\end{proposition}
\theproof
We will show that we have an equivalence
\(
\ZZ/2\times \prod_{k>0}K(\ZZ,2k)\simeq  \op{GL}_1(\H\ZZ[u,u^{-1}])\;.
\label{Z2prod}
\)
The connected cover of $\H\ZZ[u,u^{-1}]$ is given 
by $\H\ZZ[u]$ and the infinite loop space is the Dold-Kan image of the positively graded 
complex $\ZZ[u]\cong \prod_{k}\ZZ[2k]$, which is a model for $\prod_{k}K(\ZZ,2k)$. 
Since the group of units of $\ZZ$ are $\ZZ/2\cong\{-1,1\}$, we see that 
$\op{GL}_1(\H\ZZ[u,u^{-1}])$ is as claimed and delooping gives the desired equivalence. 
\endofproof

We now describe  the twists via module spectra. 
For a ring spectrum $\mathscr{R}$, let us recall the Picard $\infty$-groupoid 
${\rm Pic}^{\rm top}_{\mathscr{R}}$ from \cite{GS5}, following \cite{BN}. This is 
the infinity groupoid whose objects are invertible $\mathscr{R}$-module spectra. The corresponding 
 geometric realization decomposes in the category of spaces as
$$
\vert {\rm Pic}^{\rm top}_{\mathscr{R}} \vert \simeq  B{\rm GL}_1(\mathscr{R}) 
\times   \pi_0{\rm Pic}^{\rm top}_{\mathscr{R}} \;.
$$
For the spectrum $\H\ZZ[u,u^{-1}]$, Proposition \ref{prop twhz} then gives a canonical map
$$
\xymatrix{
K(\ZZ/2,1)\times \prod_{k>0}K(\ZZ,2k+1)
\; \ar@{^{(}->}[r]& 
\vert  {\rm Pic}^{\rm top}_{\H\ZZ[u,u^{-1}]} \vert \;,
}
$$
given by the inclusion at the identity component of the Picard space. This indeed 
allows us to twist periodic integral cohomology by any odd degree integral class
 ($\ZZ/2$-class in degree one).

\medskip
We now would like to describe the actual module spectra which exhibit the twisted theory. 
One of the most systematic ways to describe the resulting module spectra was presented 
in \cite{GS5}. There, we defined a canonical bundle of spectra over the Picard infinity 
groupoid which lives in the tangent infinity topos\footnote{$\mathscr{S}{\rm pace}$ is the category 
of compactly generated, weakly Hausdorff spaces. The only reason for this condition on our topological 
spaces is so that we have a convenient category of spaces in which to work. In particular, we have internal 
mapping spaces in this category which turns this category into an $\infty$-category.} $T(\mathscr{S}{\rm pace})$ 
and the pullbacks of this universal bundle by a map $h:X\to {\rm Pic}^{\rm top}_{\mathscr{R}}$ gave 
a bundle of spectra representing the twisted theory. Since we would like to be as concrete as possible, and 
relying on as little abstract machinery as possible, we note that in the present case this universal bundle 
will take on a relatively simple form; see the map \eqref{univ bun hquo}.


\medskip
We begin by describing a convenient category in which our bundles of spectra live. Ordinary vector 
bundles are allowed to live in topological spaces since the fiber itself is a topological space. 
However, a spectrum is a generalization of certain types of topological spaces, namely infinite 
loop spaces. An infinite loop space  is allowed to have  homotopy groups in negative degrees and, therefore,
 cannot itself be regarded as a space. As mentioned above, the convenient  category in which our constructions take place
  is the tangent infinity category of spaces,  $T(\mathscr{S}{\rm pace})$ (see Remark \ref{stinf parspec} below) . 
  Similar to the way one defines a spectrum from a prespectrum, we have the following definitions and 
  properties.\footnote{After some identifications  the definition is essentially \cite[Def. 11.2.3]{MS}, in the 
  case $G=\ast$. The corresponding $\infty$-category theoretic treatment can be found in 
  \cite[Sec. 35.5]{Joy}.}

\begin{definition}
[Parametrized (pre)spectra]
\label{Def-par-spec}
{\bf (i)} A \emph{parametrized prespectrum over $X$} is a collection of maps between spaces
$
p_n: E_n\to X
$, $n \in \ZZ$, together  with a choice of section, and which come equipped with morphisms 
$\Sigma_XE_n\to E_{n+1}$, commuting with the sections. Moreover, the 
following diagram 
$$
\xymatrix@R=1.5em{
\Sigma_XE_n\ar[rr]\ar[dr]_-{\Sigma_X(p_n)} &&  E_{n+1}\ar[dl]^-{p_{n+1}}
\\
& X& 
}
$$
is required to commute up to a choice of equivalence. 
The operation $\Sigma_X$ is the result of the pushout
$$
\xymatrix@R=1.5em{
E_n\ar[rr]^-{p_n}\ar[d]_-{p_n} &&
  X\ar[d]
\\
X\ar[rr] && \Sigma_X E_n\;.
}
$$
\item {\bf (ii)} A parametrized prespectrum $\{p_n: E_n\to X\}$ for which the adjoint maps 
$E_n\to \Omega_XE_{n+1}$ are equivalences is called 
a \emph{parametrized spectrum}. 
\end{definition}

Note that the above pushout $\Sigma_X  E_n$ is a generalization of the 
usual suspension $\Sigma E_n$, which arises by taking $X$ to be a point.

\begin{remark}
[The tangent infinity category of spaces from parametrized spectra]
\label{stinf parspec}
Forgetting about the levels of the object $\{p_n\}$, we will often denote a parametrized spectrum 
whose levels are $p_n: E_n\to X$, simply by $p: E\to X$.  Given two objects $p: E\to X$ and $p^{\prime}:
E^{\prime}\to X^{\prime}$ we have a (ordinary, unparametrized) spectrum of maps $\map(p,p^{\prime})$, 
which in particular (at level zero) gives a collection of homotopy  commutative diagrams
$$
\left\{    \vcenter{\xymatrix@R=1.3em{
E_n\ar[rr]\ar[d] && E_n^{\prime}\ar[d]\ar[d]
\\
X\ar[rr] && X^{\prime}
}}\right\}\;.
$$
The resulting structure is that of a stable $\infty$-category, which we denote by 
$T(\mathscr{S}{\rm pace})$. We denote the subcategory on those spectra 
parametrized over a fixed space $X$ as $T(\mathscr{S}{\rm pace})_X$.
\end{remark}

Now that we have a convenient category of parametrized spectra in which 
to work, we can define a bundle of spectra as follows.

\begin{definition}
[Bundle of spectra]
\label{Def-bun-spec}
A \emph{bundle of spectra} $\pi:E\to X$ over a space $X$ with fiber 
a (ring) spectrum $\mathscr{R}$
 is an object in $T(\mathscr{S}{\rm pace})$ satisfying the following properties:
\begin{enumerate}[{\bf (i)}]

\item  \emph{Fiber:}  For each $x\in X$, the pullback 
$$
\xymatrix@R=1.3em{
\pi^{-1}(x)\ar[rr] \ar[d]
&& E\ar[d]
\\
\ast\ar[rr]^x  && X
}
$$
is equivalent to $\mathscr{R}$.

\item \emph{Local trivialization:} There is a covering $\{U_{\alpha}\}$ of $X$ such that, 
for each $U_{\alpha}$, we have a Cartesian square (i.e. $(\infty,1)$-pullback square)
$$
\xymatrix@R=1.3em{
U_{\alpha}\times \mathscr{R}\ar[rr]^-{\phi_{\alpha}}_>>>>>{\ }="s" \ar[d] && E\ar[d]
\\
U_{\alpha}\ar[rr]^<<<<<{\ }="t" && X
 \ar@{=>}_{h_\a} "s"; "t"
\;.}
$$
The map $\phi_{\alpha}$ and the homotopy $h_{\alpha}$ filling the diagram constitute what we call 
a local trivialization. \footnote{Note that the map $\phi_{\alpha}$ and the homotopy $h_{\alpha}$ 
are \emph{not} enough to reconstruct the bundle. One also needs to consider pullbacks to higher-fold 
intersections and higher homotopies filling the resulting diagrams. This is in contrast to the case of an 
ordinary vector bundle, where it is enough to know the local trivializations.}
\end{enumerate}
\end{definition}
In the absence of any geometry, bundles of spectra behave more like covering spaces than like smooth 
vector bundles. The next example illustrates this point.
\begin{example}
[Bundle of spectra over the circle]
\label{exS1}
Let $Z\to S^1$ be the disconnected cover of $S^1$, splitting as the disjoint union $Z=\coprod_{k}W$, 
where $W$ is the connected cover classified by the subgroup $2\ZZ\subset \ZZ\cong \pi_1(S^1)$. This 
cover can be viewed as a $\ZZ$-subbundle of the M\"obius bundle given by restricting to integers. 
Viewing $S^1$ as the unit circle in the complex plane and removing the points $-1$ and $1$ from $S^1$, 
we get corresponding open sets $U$ and $V$, respectively, covering $S^1$. Over $U$ and $V$, we have 
equivalences
$$
\phi_{U}:Z\vert_U\simeq \ZZ \times U
\qquad \text{and} \qquad
  \phi_{V}:Z\vert_V\simeq \ZZ\times V\;,
$$
which can be chosen so that the transition functions act by multiplication by $-1$ on the fibers,
i.e., $\phi_{UV}(n,x)=(-n,x)$. The map $-1\times:\ZZ\to \ZZ$ extends to a map 
$-1\times :\H(\ZZ[u,u^{-1}])\to \H(\ZZ[u,u^{-1}])$ degreewise. Gluing by this automorphism gives the identification
$$
\left\{ \vcenter{\xymatrix{ {\mathcal Z}\ar[d]
\\
S^1
}}\right\} \simeq {\rm colim} \left\{  \vcenter{\xymatrix{
\H(\ZZ[u,u^{-1}])\times U\cap V\ar[d] \ar@<.07cm>[r]^-{-i_1}\ar@<-.07cm>[r]_-{i_2} & \H(\ZZ[u,u^{-1}])\times U \coprod \H(\ZZ[u,u^{-1}])\times V\ar[d]
\\
U\cap V \ar@<.07cm>[r]\ar@<-.07cm>[r] & U\coprod V 
}}
\right\}\;,
$$
where $i_2$
is induced by the usual inclusion into the 
second factor and the top map $-i_1$ applies the automorphism $-1$ and then includes into the first factor. 
This colimit takes place in the category $T(\mathscr{S}{\rm pace})$. As part of the data of the colimit, 
we have local trivializations
$$
\phi_{U}:{\mathcal Z}\vert_U\simeq \H(\ZZ[u,u^{-1}])\times U
\qquad \text{and} \qquad  
\phi_{V}:{\mathcal Z}\vert_V\simeq \H(\ZZ[u,u^{-}])\times V\;,
$$
turning ${\mathcal Z}$ into a corresponding bundle of spectra, with fiber $\H (\ZZ[u, u^{-1}])$, 
over $S^1$. The transition functions take the form $\phi_{UV}(x,-)=-1$, where $-1$ is the automorphism of  
the fiber $\H(\ZZ[u,u^{-1}])$ induced by multiplication by $-1$.
\end{example}

In Example \ref{exS1}, the automorphisms $\phi_{UV}(x,-)$ had degree zero, in the sense that they were 
genuine 1-morphisms and not higher simplices in the space of automorphisms ${\rm GL}_1(\H(\ZZ[u,u^{-1}))$. 
The next example, however, gives an instance where we do have higher simplices. 

\begin{example}
[Bundle of spectra over the 3-sphere]
\label{3 sphex}
Consider the 3-sphere $S^3$, equipped with the cover $\{U,V\}$ obtained by removing the north 
and south poles, respectively. The intersection $U\cap V\simeq S^2$ and, given our identification 
of the units in (the proof of) Proposition \ref{prop twhz}, we have 
$$
\pi_2({\rm GL}_1(\H\ZZ[u,u^{-1}]))\simeq \pi_2(K(\ZZ,2))\simeq \ZZ\;,
$$
with generator $u$. Consequently,  the homotopy class of a map 
$U\cap V \simeq S^2 \to {\rm GL}_1(\H\ZZ[u,u^{-1}])$ 
is represented by an integer $n$ times the generator $u$. Via the action of 
${\rm GL}_1(\H\ZZ[u,u^{-1}])$, such a representative gives rise to a map
\(\label{sphunit gl1}
nu:S^2\longrightarrow \map\big(\H\ZZ[u,u^{-1}], \;\H\ZZ[u,u^{-1}]\big)\;,
\)
and we would like to take the map as supplying the transition data for a bundle on $S^3$. 
Acting by this map  and then by the usual inclusion map $U\cap V\into U\coprod V$ into the 
second factor gives the two top arrows (that is, $\times nui_1$ and $i_2$), respectively, 
 in the following diagram
$$
\left\{ \vcenter{\xymatrix{ {\mathcal Z}\ar[d]
\\
S^3
}}\right\} \simeq {\rm hocolim} 
\left\{  \vcenter{\xymatrix{
 \H(\ZZ[u,u^{-1}])\times U\cap V
  \ar[dd] \ar@/^2pc/[rr]^-{\times nu i_1}_>{\ }="s"
  \ar@/_2pc/[rr]_-{i_2}^<<<<<<{\ \hspace{10cm}}="t"
  &&
  \H(\ZZ[u,u^{-1}])\times  U\coprod \H(\ZZ[u,u^{-1}])\times  V \ar[dd]
\\
\\
U\cap V  \ar@<.07cm>[rr]
\ar@<-.07cm>[rr]
 &&U\coprod V
{\ar@{=>}@/_.6pc/ "s"+(-34, 3); "t"}
{\ar@{=>}@/^.6pc/ "s" +(-28, 3); "t"}
{\ar@3{->} (27,0)*{}; (32,0)*{}}
}
}
\right\}\;.
$$
The fact that this diagram has nontrivial homotopies filling it, i.e., the ones
provided by the map \eqref{sphunit gl1}, is what separates it from Example \ref{exS1}.  
The homotopy class of sections of the bundle ${\mathcal Z}\to S^3$ computes
 the twisted cohomology groups. 
\end{example}

In the same way that ordinary vector bundles with $G$-structure are classified by maps to 
the classifying space $BG$, bundles of spectra with fiber $\mathscr{R}$ are classified by maps 
to $B{\rm GL}_1(\mathscr{R})$. There is a universal bundle of spectra over this space. 
In the present case (i.e. for periodic integral cohomology) it takes the following form.
The action of each factor $K(\ZZ,2k)$ on the spectrum 
$\H\ZZ[u,u^{-1}]$ gives rise to a quotient 
\footnote{Since we are in an $(\infty, 1)$-category, quotients are taken in 
the $(\infty, 1)$-sense, i.e., up to coherence homotopy.}
$\H\ZZ[u,u^{-1}]/\!/K(\ZZ,2k)$.
This leads to the following bundle  
\(\label{univ bun hquo}
\H\ZZ[u,u^{-1}]/\!/K(\ZZ,2k)\longrightarrow K(\ZZ,2k+1)\;,
\)
which we can think of as a universal bundle.
Given a map $h:X\to K(\ZZ,2k+1)$, we consider the pullback diagram
\(
\label{higher bundle}
\xymatrix@R=1.5em{
{\mathcal Z}_{h} \ar[rr]\ar[d] &&  \H\ZZ[u,u^{-1}]/\!/K(\ZZ,2k)\ar[d]
\\
X\ar[rr]^-{h} &&K(\ZZ,2k+1)
\;.}
\)
Then ${\mathcal Z}_h\to X$ is itself a bundle of spectra with fiber $\H\ZZ[u,u^{-1}]$. Indeed, 
the Pasting Lemma for pullbacks implies that we have a double pullback square
\(
\label{higher bundle}
\xymatrix@R=1.5em{
\H\ZZ[u,u^{-1}]\ar[r] \ar[d] & {\mathcal Z}_{h} \ar[rr]\ar[d] &&  \H\ZZ[u,u^{-1}]/\!/K(\ZZ,2k)\ar[d]
\\
\ast\ar[r] & X\ar[rr]^-{h} &&K(\ZZ,2k+1)
\;,
}
\)
so that $\H\ZZ[u,u^{-1}]$ is identified as the fiber. Next, suppose that 
$X$ admits a good open cover $\{U_{\alpha}\}$. Then, by the Borsuk Nerve Theorem 
(see, e.g., \cite[Theorem 3.21]{Pr}),  $X$ is homotopy equivalent to the  colimit 
over the {\v C}ech nerve of a good open cover $\{U_{\alpha}\}$. By iterating pullbacks, 
we therefore get induced commutative simplicial diagrams
\(
\label{higher bundle csimp}
\xymatrix@=2em{
\cdots \ar@<.1cm>[r]\ar[r] \ar@<-.1cm>[r] & \coprod_{\alpha\beta}\H\ZZ[u,u^{-1}]\times U_{\alpha\beta} 
 \ar@<.05cm>[r] \ar@<-.05cm>[r]\ar[d] &\coprod_{\alpha}\H\ZZ[u,u^{-1}]\times U_{\alpha} \ar[r] \ar[d] 
 & {\mathcal Z}_{h} \ar[r]\ar[d] &  \H\ZZ[u,u^{-1}]/\!/K(\ZZ,2k)\ar[d]
\\
\cdots \ar@<.1cm>[r]\ar[r] \ar@<-.1cm>[r]& 
\coprod_{\alpha\beta}U_{\alpha\beta}  \ar@<.05cm>[r] \ar@<-.05cm>[r] 
& \coprod_{\alpha}U_{\alpha} \ar[r] & X\ar[r]^-{h} &K(\ZZ,2k+1)
\;,
}
\)
where the bottom  simplicial diagram is induced by the {\v C}ech nerve. Via descent, the top simplicial 
diagram in \eqref{higher bundle csimp} is homotopy 
colimiting and this says that (up to homotopy equivalence) we can recover the total space ${\mathcal Z}_h$ 
by gluing together local trivializations via compatibility maps defined on various 
intersections.\footnote{The tangent $\infty$-category of spaces is an example of an $\infty$-topos and such
 infinity categories are characterized axiomatically 
via the Giraud-Rezk-Lurie axioms \cite[Sec. 6.1.5]{Lur}. One of these 
axioms is that of descent, which asserts that whenever we have a diagram of the above form with the bottom 
simplicial diagram being colimiting, and all squares  being Cartesian, then the top simplicial 
diagram is also colimiting.} This association gives the following correspondence.
\begin{proposition}
[Characterization of twisted periodic integral cohomology]
There is a bijective correspondence between homotopy classes of maps $h:X\to K(\ZZ,2k+1)$ and equivalence 
classes of bundles of spectra with fiber $\H\ZZ[u,u^{-1}]$, which admit a $K(\ZZ,2k)$-structure, i.e., 
a reduction of the structure $\infty$-group from ${\rm GL}_1(\H\ZZ[u,u^{-1}])$ to $K(\ZZ, 2k)$.
\end{proposition}

\begin{example}
[Classifying map for bundles of spectra over $S^3$.]
In Example \ref{3 sphex} the transition data specified by the map $\times nu: S^2\simeq U\cap V\to K(\ZZ,2)
\into {\rm GL}_1(\H\ZZ[u,u^{-1}])$ corresponds to a map 
$h:S^3\to K(\ZZ,3)\into B{\rm GL}_1(\H\ZZ[u,u^{-1}])$ by the loop-suspension adjunction. 
This map is the classifying map of the bundle constructed in that example.
\end{example}

As we stated in Remark \ref{stinf parspec}, the sections of the map $p:{\mathcal Z}_h\to X$ form a spectrum.
 Given that, locally, ${\mathcal Z}_h$ trivializes as $\H\ZZ[u,u^{-1}]\times U_{\alpha}$ when $\{U_{\alpha}\}$
  is a good open cover of a space $X$, we can calculate the spectrum via the local data as the  limit of spectra
\(
\Gamma(X;{\mathcal Z}_h)={\rm lim}\left\{ \xymatrix{
\cdots & \ar@<.1cm>[l]\ar[l] \ar@<-.1cm>[l]  \coprod_{\alpha\beta}\H\ZZ[u,u^{-1}]& \ar@<.05cm>[l] \ar@<-.05cm>[l] \coprod_{\alpha}\H\ZZ[u,u^{-1}]
}\right\}\;,
\)
where again the simplicial homotopy commutative diagram is determined by the transition functions 
and higher transition data. In practice, this can aid in calculation; however, it is more useful to 
develop some of the basic properties of the spectrum of sections. Indeed, we will do this in Section \ref{Sec-prop}. 

\medskip
We finish our current discussion by defining the  underlying twisted cohomology groups for twisted periodic 
$\ZZ$-cohomology.  Notice that, since the fibers $\H\ZZ[u,u^{-1}]$ are 2-periodic, in the sense that 
$\Sigma^2\H\ZZ[u,u^{-1}]\simeq \H\ZZ[u,u^{-1}]$, and the action by $K(\ZZ,2k)$ commutes with this shift, the 
sections of ${\mathcal Z}_h$ are also 2-periodic. This leads us to the following definition for the reduced 
cohomology.

\begin{definition}
[Twisted periodic integral cohomology]
Let $h:X\to K(\ZZ,2k+1)$ be be a twist for periodic integral cohomology. We define the 
\emph{$h$-twisted integral cohomology} as the $\ZZ/2$-graded group
$$
\widetilde{H}^*(X;h):=\pi_*\Gamma(X,{\mathcal Z}_h)\;.
$$
We will refer to the degree of a class as either \emph{even} or \emph{odd}, 
corresponding to the identity and nonidentity elements in $\ZZ/2$, respectively.
\end{definition}

\subsection{Properties of twisted periodic integral cohomology}
\label{Sec-prop}

In this section, we state some of the basic properties of twisted periodic cohomology, 
which we generalize to twisted periodic smooth Deligne cohomology in Section 
\ref{Sec-cech}. The following proposition holds more generally for 
any twisted cohomology theory, but we will only state this in the present case for the reduced
theory $ \widetilde{H}^*(X;h)$.

\begin{proposition}
[Properties of twisted periodic integral cohomology]
\label{prop-tpic}
Let $X$ be a space and fix a twist as a map $h:X\to K(\ZZ,2k+1)$. Consider the category of such 
pairs $(X,h)$, with morphisms $f:(X,h)\to (Y, \ell)$ given by maps $f:X\to Y$ such that 
$[f^*\ell]=[h]$. The assignment $(X,h)\mapsto \widetilde{H}^*(X;h)$ satisfies the following 
properties:
\begin{enumerate}[{\bf (i)}]
\item $\widetilde{H}^*(M;h)$ is functorial with respect to the maps $f:(X,h)\to (Y,\ell)$.

\item The functor $\widetilde{H}^*(-;h)$ satisfies the Eilenberg-Steenrod axioms for a 
reduced generalized cohomology cohomology theory 
(i.e., modulo the dimension axiom). In particular, we have a Mayer-Vietoris sequence
\begin{equation*}
\xymatrix{
\widetilde{H}^{\rm ev }(M;h)\ar[r] & \widetilde{H}^{\rm ev}(U;h)\oplus\widetilde{H}^{\rm ev}(V;h)\ar[r]& 
\widetilde{H}^{\rm ev}(U\cap V ;h) \ar[d]^-{\partial}
\\ 
\widetilde{H}^{\rm odd}(U\cap V;h)\ar[u]^-{\partial}
&\ar[l] \widetilde{H}^{\rm odd}(U;h)\oplus\widetilde{H}^{\rm odd}(V;h)& 
\ar[l] \widetilde{H}^{\rm odd}(M;h)
 }
\end{equation*}
where $\partial$ is the connecting homomorphism, and the sequence is exact at each entry. 

\item For $h:X\to K(\ZZ,2k+1)$ a trivial twist (i.e. $h\simeq \ast$) we have an isomorphism
$$
\widetilde{H}^*(X;h)\cong \widetilde{H}^*(X;\ZZ[u,u^{-1}])\;.
$$
\end{enumerate}
\end{proposition}

\theproof
{\bf (i)} Given a map $f:X\to Y$ satisfying the desired compatibility, we have an induced double pullback diagram 
\(
\label{higher bundle czh}
\xymatrix@R=1.5em{
 {\mathcal Z}_{f^*\ell}\ar[r]\ar[d] & {\mathcal Z}_{\ell} \ar[rr]\ar[d] &&  
 \H\ZZ[u,u^{-1}]/\!/K(\ZZ,2k)\ar[d]
\\
X \ar[r]^f & Y\ar[rr]^-{\ell} &&K(\ZZ,2k+1)
\;,
}
\)
which gives an identification ${\mathcal Z}_{f^*\ell}\simeq {\mathcal Z}_h$.  Consequently, 
we have an induced morphism of sections 
$f^*:\Gamma(X;{\mathcal Z}_\ell)\to \Gamma(Y;{\mathcal Z}_{f^*\ell}) 
\overset{\simeq}{\to} \Gamma(Y, {\mathcal Z}_h)$. Passing to homotopy groups yields a map 
$f^*:\widetilde{H}^*(Y; \ell)\to \widetilde{H}^*(X; h)$. It is clear that this assignment takes compositions
 of morphisms of pairs to compositions of group homomorphisms.

\medskip
\noindent {\bf (ii)} We now verify the generalized Eilenberg-Steenrod axioms. 
\vspace{0.1cm}

\noindent {\it Homotopy invariance}. It follows from the universal property of the pullback that when two maps 
$f:(X,h)\to (Y, \ell)$ and $g:(X,h)\to (Y, \ell)$ are homotopic, \footnote{Note that homotopy is a relation between 
morphisms of \emph{pairs} $(X,h)$ and must respect the maps to the space of twists.} the induced map on sections 
$f^*:\Gamma(X;{\mathcal Z}_\ell)\to \Gamma(Y;{\mathcal Z}_{f^*\ell}) 
\overset{\simeq}{\to} \Gamma(Y, {\mathcal Z}_h)$ and 
$g^*:\Gamma(X;{\mathcal Z}_\ell)\to \Gamma(Y;{\mathcal Z}_{g^*\ell}) 
\overset{\simeq}{\to} \Gamma(Y, {\mathcal Z}_h)$ are homotopic and, therefore,
 induce isomorphic maps at the level of cohomology.

\vspace{0.1cm}
\noindent {\it Additivity}. 
Let $X=\coprod_{\alpha}X_{\alpha}$, then a map $h:X\to K(2k+1)$ is equivalently a collection of maps 
$h_{\alpha}:X_{\alpha}\to K(\ZZ,2k+1)$. Hence, the spectrum of sections $\Gamma(X,{\mathcal Z}_{h})$ 
splits as a product $\prod_{\alpha}\Gamma(X_{\alpha},{\mathcal Z}_{h_{\alpha}})$. 
Since taking homotopy groups commutes with products, we have an isomorphism
$$
\widetilde{H}^*(X,h)\cong \prod_{\alpha}\widetilde{H}^*(X_{\alpha},h_{\alpha})\;.
$$
\noindent {\it Exactness}. Let $i:(A,i^*h)\into (X,h)$ be an inclusion and consider the cofiber sequence 
$(A,i^*h)\into (X,h)\to {\rm cone}(i,\tilde{h})$. The functor $\Gamma(-;{\mathcal Z})$ sends this homotopy
 cofiber sequences to homotopy fiber sequences and, therefore, we have a fiber sequence 
$
\Gamma({\rm cone}(i);{\mathcal Z}_{\tilde{h}})\to
\Gamma(X;{\mathcal Z}_h) 
\to
 \Gamma(A;{\mathcal Z}_{i^*h})
 $.
The associated long exact sequence for the cofiber sequence 
$$
\xymatrix{
(U\cap V,r_{{}_{UV}}^*h) \; \ar@{^{(}->}[r] &
 (U,r_{{}_U}^*h)\vee_{K(\ZZ,2k+1)}(V,r_{{}_V}^*h)\ar[r] & (U\cup V,h)}\;,
$$ 
with $r_{{}_W}$ the restriction to the appropriate open set $W$, gives the Mayer-Vietoris sequence. 

\noindent {\bf (iii)} Finally, if the twist $h:X\to K(\ZZ,2k+1)$ is trivial then, up to homotopy,
 $h$ factors through the point inclusion $\ast\into K(\ZZ,2k+1)$ induced by the zero map
$0 \to \ZZ$. Fixing such a homotopy then taking iterated  pullbacks and using the Pasting 
Lemma gives  pullback squares 
\(
\label{dpull triv}
\xymatrix@R=1.5em{
 {\mathcal Z}_{h}\simeq \H\ZZ[u,u^{-1}]\times X\ar[r]\ar[d] &  \H\ZZ[u,u^{-1}] \ar[rr]\ar[d] 
 &&  \H\ZZ[u,u^{-1}]/\!/K(\ZZ,2k)\ar[d]
\\
X \ar@/_1pc/[rrr]_{h}\ar[r]  & \ast \ar[rr] &&K(\ZZ,2k+1)
\;.}
\)
The bundle equivalence ${\mathcal Z}_h\simeq \H\ZZ[u,u^{-1}]\times X$ then induces an 
equivalence at the level of global sections, and hence an isomorphism at the level of 
corresponding reduced theories.
\endofproof

\begin{remark} 
[Reduced vs. unreduced]
Note that we can pass from the reduced theory $\widetilde{H}^*(X;h)$ to the unreduced 
theory of a triple $(X,A,h)$ with $A\subset X$ as usual, by setting 
$$
H^*(X,A,h)\cong \widetilde{H}^*(X/A;h)\oplus \ZZ[u,u^{-1}]\;.
$$ 
\end{remark}

\section{Periodic smooth Deligne cohomology}
\label{Sec Per Del}

In this section, we introduce the notion of periodic Deligne cohomology. This will set the stage for the next section, 
where we identify the twists of this theory. 

\subsection{Construction as a cohomology theory}
\label{Sec-ccoh}

Just as the Deligne complex is indexed by an integer $n\in \ZZ$, here we have complexes indexed by elements in $\ZZ/2$, 
which we will call either even or odd, depending on the parity. 
We let $\underline{\ZZ}$ denote the locally constant sheaf of $\ZZ$-valued functions.

\begin{definition}
[Even and odd Deligne complexes]
\label{Devodd}
For $\op{ev}$, $\op{odd}\in \ZZ/2$, corresponding to the identity and nonidentity components, respectively, we have the 
two complexes
$$
{\mathcal D}(\op{ev}):=\big(\underbrace{\hdots \longrightarrow  \underline{\ZZ}\oplus \prod_{k}\Omega^{2k+1}\longrightarrow  \prod_{k}\Omega^{2k}\longrightarrow \underline{\ZZ}\oplus \prod_{k}\Omega^{2k+1}}_{{\rm deg} \geq 0}\longrightarrow \underbrace{0 \longrightarrow    \vphantom{  \prod_{k}\Omega^{2k+1} }\underline{\ZZ} \longrightarrow \hdots}_{{\rm deg}< 0}
\big)
$$
and 
$$
{\mathcal D}(\op{odd}):=\big(\underbrace{\hdots \longrightarrow  \prod_{k}\Omega^{2k}\longrightarrow   \underline{\ZZ}\oplus \prod_{k}\Omega^{2k+1}\longrightarrow \prod_{k}\Omega^{2k}}_{{\rm deg} \geq 0}\longrightarrow \underbrace{0 \longrightarrow   \vphantom{  \prod_{k}\Omega^{2k+1} }\underline{\ZZ} \longrightarrow \hdots}_{{\rm deg}< 0}
\big)\;,
$$
where the $\underline{\ZZ}$'s sit in even degrees in the first complex and in odd degrees in the second. In positive degrees, 
the differential in both complexes is the usual exterior derivative term-wise and on the copies of $\underline{\ZZ}$ it is given by the inclusion map $\underline{\ZZ}\into \Omega^0\into \prod_{k}\Omega^{2k}$. In negative degrees the differential is trivial.
\end{definition}

The complexes ${\mathcal D}(\op{ev})$ and ${\mathcal D}(\op{odd})$ are sheaves of chain complexes on 
the category of all smooth manifolds, topologized as a site via good open covers. Alternatively, 
both complexes can be regarded as sheaves of chain complexes 
on any fixed manifold $M$ simply by evaluating on the open subsets of $M$. 
This is the familiar setting in which ordinary smooth Deligne cohomology takes 
place (e.g. \cite{Bry}). We have the following natural definition.
 
\begin{definition}
[Periodic Deligne cohomology]
\label{def pdc}
We define \emph{$\ZZ/2$-graded periodic Deligne cohomology groups}
 of a smooth manifold $M$ as the sheaf  hypercohomology groups
 \footnote{The cohomological degrees on both right hand sides is $0$, due to 
 the shift in the complexes in Definition \ref{Devodd}, in analogy to the usual Deligne 
 case, i.e., expression \eqref{Delcplx}.}
$$
\widehat{H}^{\op{ev}}(M;\ZZ[u,u^{-1}]):=H^0(M;{\mathcal D}(\op{ev})) 
\qquad \text{and} \qquad
 \widehat{H}^{\op{odd}}(M;\ZZ[u,u^{-1}]):=H^0(M;{\mathcal D}(\op{odd}))\;.
$$
\end{definition}

The following shows that periodic Deligne cohomology can be calculated easily from the ordinary Deligne cohomology 
groups of a manifold.

\begin{proposition}
[Calculating periodic Deligne cohomology groups]
\label{spl pdel}
Let $M$ be a smooth manifold. There are natural isomorphisms
$$
\widehat{H}^{\op{ev}}(M;\ZZ[u,u^{-1}])\cong \bigoplus_{k}\widehat{H}^{2k}(M;\ZZ)
\qquad \text{and} \qquad
\widehat{H}^{\op{odd}}(M;\ZZ[u,u^{-1}])\cong \bigoplus_{k}\widehat{H}^{2k+1}(M;\ZZ)\;.
$$
\end{proposition}
\theproof
We will prove the claim for $\widehat{H}^{\op{ev}}(M;\ZZ[u,u^{-1}])$. The case for 
$\widehat{H}^{\op{odd}}(M;\ZZ[u,u^{-1}])$ is proved similarly. To this end, we 
organize the sheaf of chain complexes ${\mathcal D}(\op{ev})$ as follows
$$
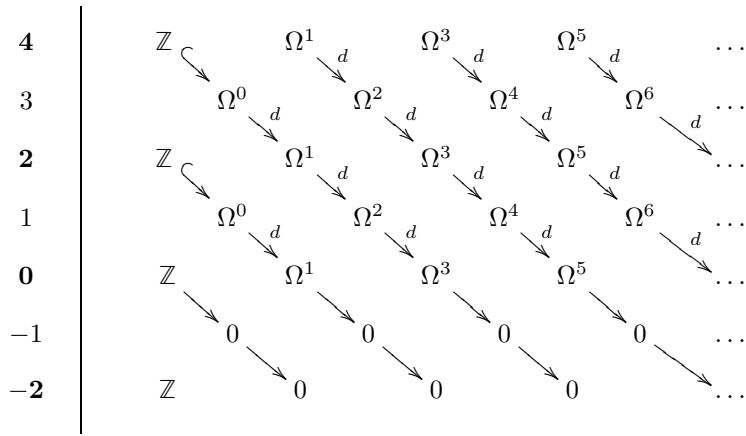
\xymatrix@=8pt{
& \ar@{-}[dddddddd]  & & & & & & & & & & 
\\
{\bf 4} & & & \ZZ \; \ar@{^{(}->}[dr] & & \Omega^1\ar[dr]^-{d} & & \Omega^3 \ar[dr]^-{d} & & \Omega^5 \ar[dr]^-{d}  & & \; \; \;\hdots 
\\
3 & & &  & \Omega^0 \ar[dr]^-{d}& & \Omega^2  \ar[dr]^-{d}& & \Omega^4  \ar[dr]^-{d}& & \Omega^6 \ar[dr]^-{d} &  \; \; \; \hdots
 \\
{\bf 2} & & & \ZZ \; \ar@{^{(}->}[dr]& & \Omega^1 \ar[dr]^-{d} & & \Omega^3 \ar[dr]^-{d} & & \Omega^5  \ar[dr]^-{d}& & \; \; \; \hdots 
 \\
1 & & &  & \Omega^0 \ar[dr]^-{d}  & & \Omega^2 \ar[dr]^-{d} & & \Omega^4 \ar[dr]^-{d} & & \Omega^6 \ar[dr]^-{d}  & \; \; \; \hdots
  \\
{\bf 0} && &   \ZZ \ar[dr] & & \Omega^1 \ar[dr]& & \Omega^3 \ar[dr] & & \Omega^5\ar[dr] & &
\; \; \; \hdots 
\\
-1 & & & & 0\ar[dr] & & 0 \ar[dr]& & 0\ar[dr] & & 0\ar[dr] & 
\; \; \; \hdots
\\
{\bf -2} & & & \ZZ && 0 & & 0 & & 0 & & 
\; \; \;\hdots
\\
& & & & & & & & & & & 
}
$$
where the numbers on the vertical axis index the degree of the complex. 
The diagonal arrows represent  the differential on each component of 
the product taken over a given  row. 
The diagonal complexes are easily seen to be the usual Deligne complex and, therefore,
 we have a splitting 
\(\label{spl del e}
{\mathcal D}(\op{ev})\cong \prod_{k}{\mathcal D}(2k)\oplus \prod_{k}\ZZ[-2k]\;,
\) 
where the second summand comes from the negative degrees of the complex. 
The latter do not contribute to the hypercohomology of the complex, as the {\v C}ech resolution of the complex necessarily vanishes in negative degrees. Thus, the hypercohomology groups split as desired.
\endofproof

From Proposition \ref{spl pdel}, it follows immediately that the periodic Deligne cohomology 
groups fit into a differential cohomology diamond diagram and into exact sequences similar 
to those for ordinary Deligne cohomology, as an instance of differential integral cohomology \cite{SSu}. 

\begin{proposition}[Periodic  Deligne cohomology diamond]  
\label{perd diam diag}
 \label{periodic diamond} We have the exact diamond diagram 
\(\hspace{-.1cm}\label{diamond2}
\xymatrix@C=-1.3em {
 &\Omega^{\op{odd}}(M)/{\rm im}(d) \ar[rd]^{a}\ar[rr]^-{d} & & 
\Omega^{\op{ev}}_{\rm cl}(M) \ar[rd] &  
\\
H^{\op{odd}}(M;\RR[u,u^{-1}])\ar[ru]\ar[rd] & & 
\widehat{H}^{\op{ev}}(M;\ZZ[u,u^{-1}]) \ar[rd]^{I} \ar[ru]^{R}& &  
H^{\op{ev}}(M;\RR[u,u^{-1}]) 
\\
&H^{\op{odd}}(M;\RR[u,u^{-1}]/\ZZ[u,u^{-1}]) \ar[ru]\ar[rr]^-{\beta} & & 
H^{\op{ev}}(M;\ZZ[u,u^{-1}])  \ar[ru]^-{j} &
}
\)
for the even Deligne complex and a similar diamond for the odd one, given by switching $\op{ev}$
 and $\op{odd}$. Here $\Omega^{\op{odd}}(M)$ and $\Omega^{\op{ev}}(M)$ are the 
 groups of differential forms of odd and even degrees, respectively. For example, an 
 element $\omega\in \Omega^{\op{ev}}(M)$ is a formal combination
$$
\omega=\omega_0+\omega_2+\omega_4+\cdots\;,
$$
with $\omega_{2i}$ a differential form of degree $2i$.
\end{proposition}

\begin{remark}
[Extension of the diamond to a long exact sequence]
\label{two otlex}
One of the diagonals in the diamond diagram in Proposition \ref{perd diam diag} can be extended 
to a long exact sequence. Depending on the parity, the relevant 
segments of this long exact sequence are given by
\begin{align*}
&H^{\op{ev}}(M;\ZZ[u,u^{-1}])\longrightarrow \Omega^{\op{ev}}(M)/{\rm im}(d)\longrightarrow 
\widehat{H}^{\op{odd}}(M;\ZZ[u,u^{-1}]) \longrightarrow H^{\op{odd}}(M;\ZZ[u,u^{-1}])\longrightarrow 0\;,
\\
&H^{\op{odd}}(M;\ZZ[u,u^{-1}])\longrightarrow \Omega^{\op{odd}}(M)/{\rm im}(d)\longrightarrow \widehat{H}^{\op{ev}}(M;\ZZ[u,u^{-1}]) \longrightarrow H^{\op{ev}}(M;\ZZ[u,u^{-1}])\longrightarrow 0\;.
\end{align*}
The map into the quotient $\Omega^{\op{ev}}(M)/{\rm im}(d)$ takes a periodic integral class and maps it to the class of its corresponding de Rham representative (i.e. a form with integral periods). Note also that the map $R:\widehat{H}^{\op{ev}}(M;\ZZ[u,u^{-1}])\to \Omega^{\rm ev}_{\rm cl}(M)$ is not surjective; 
its image is the subgroup of closed forms with integral periods.
\end{remark}

\subsection{Ring structure and examples}
\label{Sec-ring}

Eventually, we would like to consider the twists of this theory and, to do this, we need 
a ring structure on this periodic Deligne complex. Recall that for ordinary Deligne 
cohomology, the Deligne-Beilinson cup product gives a collection of morphisms of 
sheaves of chain complexes
\cite{Del}\cite{Be} (see also \cite{FSS1}\cite{FSS2})
\(
\label{cup-DB}
\cup_{\rm DB}:{\mathcal D}(n)\otimes {\mathcal D}(m)\longrightarrow {\mathcal D}(n+m)\;.
\)
At the level of local sections, it is defined by the formula
$$
\alpha\cup_{\rm DB}\beta=\left\{
\begin{array}{cl}
\alpha\beta, & {\rm deg}(\alpha)=n
\\
\alpha\wedge d\beta, & {\rm deg}(\beta)=0, {\rm deg}(\alpha)\neq n
\\
0, &  {\rm otherwise}.
\end{array}\right.
$$
Since the even periodic Deligne complex split as the product \eqref{spl del e} (and similarly for the odd), 
there are multiplication maps
\begin{equation}
\label{mult pdel del}
\cup_{\rm DB} :
\left\{
\begin{array}{@{}l}
{\mathcal D}(\op{ev})\otimes {\mathcal D}(\op{ev})\longrightarrow{\mathcal D}(\op{ev})\;,
 \\
{\mathcal D}(\op{ev})\otimes {\mathcal D}(\op{odd})\longrightarrow {\mathcal D}(\op{odd})\;,
 \\
{\mathcal D}(\op{odd})\otimes {\mathcal D}(\op{odd})\longrightarrow {\mathcal D}(\op{ev})\;,
\end{array}
\right.
\end{equation}
induced by the cup product $\cup_{\rm DB}$ from \eqref{cup-DB} in positive degrees and the multiplication of integers 
in negative degrees. 
It is immediate that these maps descend to a graded commutative cup product which is compatible with the 
Deligne-Beilinson cup product term-wise. We summarize these observations as follows.

\begin{proposition}
[Superalgebra structure on periodic Deligne cohomology] 
With the multiplication maps \eqref{mult pdel del} 
induced by the Deligne-Beilinson cup product, the complex ${\mathcal D}(\op{ev})\oplus {\mathcal D}(\op{odd})$ 
admits the structure of a sheaf of differentially graded superalgebras. At the level of hypercohomology, it gives 
$$
\widehat{H}^{\op{ev}/\op{odd}}(M;\ZZ[u,u^{-1}]):=\widehat{H}^{\op{ev}}(M;\ZZ[u,u^{-1}])\oplus \widehat{H}^{\op{odd}}(M;\ZZ[u,u^{-1}])
$$ 
the structure of a commutative superalgebra. Moreover, we have commutative diagrams
$$
\xymatrix{
\widehat{H}^{\op{ev}/\op{odd}}(M;\ZZ[u,u^{-1}])\otimes \widehat{H}^{\op{ev}/\op{odd}}(M;\ZZ[u,u^{-1}])
\ar[rr]^-{\cup_{\rm DB}}\ar[d] 
&& \widehat{H}^{\op{ev}/\op{odd}}(M;\ZZ[u,u^{-1}])\ar[d]
\\
H^{\op{ev}/\op{odd}}(M;\ZZ[u,u^{-1}])\otimes H^{\op{ev}/\op{odd}}(M;\ZZ[u,u^{-1}])
\ar[rr]^-{\cup} 
&& H^{\op{ev}/\op{odd}}(M;\ZZ[u,u^{-1}])
}
$$
and 
$$
\xymatrix{
\widehat{H}^{\op{ev}/\op{odd}}(M;\ZZ[u,u^{-1}])\otimes \widehat{H}^{\op{ev}/\op{odd}}(M;\ZZ[u,u^{-1}])
\ar[rr]^-{\cup_{\rm DB}}\ar[d] && 
\widehat{H}^{\op{ev}/\op{odd}}(M;\ZZ[u,u^{-1}])\ar[d]
\\
\Omega_{\rm cl}^{\op{ev}/\op{odd}}(M)\otimes \Omega_{\rm cl}^{\op{ev}/\op{odd}}(M)\ar[rr]^-{\wedge} &&
 \Omega_{\rm cl}^{\op{ev}/\op{odd}}(M)\;,
}
$$
where $H^{\op{ev}/\op{odd}}(M;\ZZ[u,u^{-1}])$ is periodic integral cohomology, endowed with the superalgebra 
structure inherited from the cup product, and $\Omega^{\op{ev}/\op{odd}}(M)$ is the superalgebra of graded differential forms.
\end{proposition}

We now illustrate this with the case of spheres. 

\begin{example}
[Periodic Deligne cohomology of even spheres]
\label{Ex-PD-S2k}
Let $S^{2k}$ be the smooth even-dimensional sphere. 
The underlying periodic integral 
cohomology is readily computed as
$$
H^{\op{ev}}(S^{2k};\ZZ[u,u^{-1}])\cong \ZZ \oplus \ZZ
\qquad 
\text{and}
\qquad
 H^{\op{odd}}(S^{2k};\ZZ[u,u^{-1}])\cong 0\;.
$$
Given the two long exact sequences in Remark \ref{two otlex}, we easily compute 
$$ 
\widehat{H}^{\op{ev}}(S^{2k};\ZZ[u,u^{-1}])\cong  \Omega^{\op{odd}}(S^{2k})/{\rm im}(d)\oplus \ZZ\oplus \ZZ
\quad
\text{and}
\quad
\widehat{H}^{\op{odd}}(S^{2k};\ZZ[u,u^{-1}])\cong 
 \Omega^{\op{ev}}(S^{2k})/\Omega^{\op{ev}}_{\rm cl,\ZZ}(S^{2k})\;,
$$
where $\Omega^{\op{ev}}_{\rm cl,\ZZ}(S^{2k})$ is the subgroup of closed even forms with integral periods 
(i.e., each component of an element $\omega=\omega_0+\omega_2+\cdots$ has integral periods).
\end{example}

\begin{example}
[Periodic Deligne cohomology of odd spheres]
Similarly, we calculate for odd spheres using the same two sequences above, to get in this case
$$ \widehat{H}^{\op{odd}}(S^{2k+1};\ZZ[u,u^{-1}])\cong 
 \Omega^{\op{ev}}(S^{2k+1})/{\rm im}(d)\oplus \ZZ
$$
and
$$
 \widehat{H}^{\op{ev}}(S^{2k+1};\ZZ[u,u^{-1}])\cong 
  \Omega^{\op{odd}}(S^{2k+1})/\Omega^{\op{odd}}_{\rm cl,\ZZ}(S^{2k+1})\oplus \ZZ\;,
  $$
  where one of the $\ZZ$ factors has moved,  in comparison to the case of 
  even spheres, due to parity reasons. 
\end{example}

\section{Twisted periodic smooth Deligne cohomology} 
\label{Sec Higher}

In this section we turn to twisting periodic Deligne cohomology constructed above. 
Just as twisted periodic integral cohomology in Section \ref{Sec-bunSp} 
takes the form of a bundle of spectra over a parametrizing space, here we will have a smooth bundle 
of spectra, parametrized over a smooth manifold $M$. In the smooth setting, our starting point is 
no longer the category of spaces and its tangent infinity category $T(\mathscr{S}{\rm pace})$, but 
rather the category of smooth stacks $\sh_{\infty}(\mathscr{M}{\rm f})$ and its tangent 
infinity category $T(\sh_{\infty}(\mathscr{M}{\rm f}))$. 


\subsection{The parametrized spectrum and gerbes via smooth stacks}

Let $M$ be a smooth manifold and consider the site of open subsets $\mathcal{O}{\rm pen}(M)$, topologized 
via the good open covers $\{U_{\alpha}\to M\}$. Smooth stacks on $M$ are similar to smooth sheaves, but instead 
of assigning a set of elements to an object $U\in \mathcal{O}{\rm pen}(M)$, we assign a space (usually modeled 
combinatorially by a simplicial set). The sheaf gluing condition is replaced by a weaker condition, where we only 
require gluing up to equivalence.
To ease the transition, we start with the following. 

\begin{example}
[Gerbe with $\underline{U(1)}$-band]
\label{bry gerbe}
Consider a local homeomorphisms $U\into M$, with $U$ an open subset of $M$. To each such map, we assign 
the groupoid of line bundles with connection. This defines a Dixmier-Douady sheaf of groupoids  $\mathscr{G}$ 
on $M$. This gerbe has $\underline{U(1)}$-band, and a connective structure on this gerbe is given by a choice 
of $\Omega^1$-torsor (satisfying some properties). Using the affine structure on the space of connections, 
we see that the sheaf of connections defines 
such a torsor, so that $\mathscr{G}$ admits a connective structure. A curving of this structure is an 
assignment to each line bundle with connection $(L,\nabla)$, a two-form $K(\nabla)$ to be thought of 
as the curvature. The passage from Brylinski's gerbe \cite{Bry} to higher smooth stacks is essentially obtained 
simply by taking the nerve of the sheaf of groupoids $\mathscr{G}$, yielding an $\infty$-groupoid (which 
is a combinatorial model for a space), although one needs to be  careful in keeping track of the connective 
data and curving (see Example \ref{bry to urs} below).
\end{example}

There is a large $\infty$-category of smooth stacks $\sh_{\infty}(\mathscr{M}{\rm f})$ which does not 
depend on a choice of underlying smooth manifold. The site for this $\infty$-category is the site of all smooth 
manifolds $\mathscr{M}{\rm f}$, topologized via good open covers. 
Any object in $\sh_{\infty}(\mathscr{M}{\rm f})$ can be restricted to a single manifold by simply considering 
its value on open subsets $U\into M$. One of the benefits of working in this larger $\infty$-category is that 
one can define moduli stacks ${\bf X}$ which represent objects of interest over $M$ via maps $M\to {\bf X}$. 
For example, the moduli stack of higher gerbes with connection $\BB^nU(1)_{\nabla}$ was studied
 in \cite{Cech}\cite{SSS3}\cite{Urs}\cite{FSS2}\cite{E8}. 
 One way to present this stack is by applying the Dold-Kan functor to the sheaf of chain complexes
 $$
\BB^{n}U(1)_{\nabla}= {\rm DK} \big(\xymatrix{
 \hdots \ar[r] & 0\ar[r] & \underline{U(1)}\ar[r]^-{d\log} & \Omega^1\ar[r]^-{d} 
 & \Omega^2\ar[r]^-{d} & \hdots\ar[r] & \Omega^{n} 
 }\big)\;,
 $$
 where the sheaf $\underline{U(1)}:=C^{\infty}(-;U(1))$ sits in degree $n$. The sheaf in the argument of 
 $\op{DK}$ is quasi-isomorphic (via the exponential map) to the smooth Deligne complex ${\mathcal D}(n)$. 
 The Dold-Kan functor sends quasi-isomorphisms to weak equivalences and (since we are working up to 
 equivalence) this justifies the uniform notation $\BB^nU(1)_{\nabla}$ for both of the resulting stacks 
 (i.e. upon applying $\op{DK}$ to either complex). 
 The stack $\BB^nU(1)_{\nabla}$ sits in a Cartesian square
\(\label{stack gpull}
\xymatrix@R=1.5em{
\BB^{n}U(1)_{\nabla}\ar[rr]\ar[d] && \Omega^{n+1}_{\rm cl}\ar[d]
\\
\BB^{n+1}\ZZ\ar[rr] && \BB^{2k+1}\RR\simeq \Omega_{\rm cl}^{\leq n+1}\;,
}
\)
where $\Omega^{\leq n+1}_{\rm cl}$ is the stack presented by the sheaf of chain complexes
\(\label{de rham leq}
 \big(\xymatrix{
 \hdots \ar[r] & 0\ar[r] & \Omega^0\ar[r] & \Omega^1\ar[r]^-{d} & \Omega^2\ar[r]^-{d} &
  \hdots\ar[r] & \Omega^{n+1}_{\rm cl} 
 }\big)\;.
\)
In \cite{Cech}, it was shown that the homotopy classes of maps $M\to \BB^nU(1)_{\nabla}$ 
is in bijective correspondence with the Deligne cohomology group $\widehat{H}^n(M;\ZZ)$. 

\begin{example}
[Stack of 2-bundles with connections/gerbes with connections]
\label{bry to urs}
The smooth stack $\BB^2U(1)_{\nabla}$ can be presented via the Dold-Kan correspondence by 
the sheaf of chain complexes
$$
\BB^2U(1)_{\nabla}=L\circ \op{DK}\big(\underline{U(1)}\xrightarrow{d\log}  
\Omega^1\xrightarrow{\;d\;} \Omega^2\big)\;,
$$
where  $L$ is the stackification functor. 
\footnote{This is a functor which turns a prestack into a stack, analogously to 
the way a sheafification functor turns a presheaf into a sheaf. See 
\cite[Sec. 6.5.3]{Lur} for details.}
Let $\phi:\RR^n\to M$ 
be a local chart. For a convex open subset $U\subset \RR^n$, this stack can be evaluated on 
 the corresponding open subset $V=\phi(U)$ via
$$
\map(V,\BB^2U(1)_{\nabla})\simeq \op{DK}
\big(C^\infty(V, {U(1)})\xrightarrow{d\log}  
\Omega^1 (V)\xrightarrow{\;d\;} \Omega^2(V)\big)\;.
$$
More generally, descent for the stack $\BB^2U(1)_{\nabla}$ implies that, for any choice of good open 
cover $\{U_{\alpha}\}$ of $M$, the space of maps $\map(M,\BB^2U(1)_{\nabla})$ can be identified
 by replacing $M$ with the {\v C}ech nerve $\check{C}(\{U_{\alpha}\})$ of $\{ U_\a\}$ and considering
  instead the space of maps
$$
\map\Big(\check{C}(\{U_{\alpha}\}), \op{DK}\big(
\underline{U(1)}\xrightarrow{d\log}  
\Omega^1\xrightarrow{\;d\;} \Omega^2\big)\Big)\;.
$$
By the basic properties of the Dold-Kan correspondence we have an isomorphism 
\footnote{The shift in degree occurs because on the left we consider the complex 
$\underline{U(1)}\xrightarrow{d\log}  
\Omega^1\xrightarrow{\;d\;} \Omega^2$ 
as being shifted up two degrees relative to the complex appearing on the right.}
$$
\pi_0\Big(\map\Big(\check{C}(\{U_{\alpha}\}), \op{DK}\big(
\underline{U(1)}\xrightarrow{d\log}  
\Omega^1\xrightarrow{\;d\;} \Omega^2\big)\Big)\Big)\cong H^2(M;
\underline{U(1)}\xrightarrow{d\log}  
\Omega^1\xrightarrow{\;d\;} \Omega^2\big)\;.
$$
By \cite[Theorem 5.3.11]{Bry}, the elements on the right parametrize the homotopy classes of the gerbes 
with connective structure and curving considered in Example \ref{bry gerbe}.
\end{example}

\medskip
The definition of parametrized spectra in the smooth setting is a direct extension of 
Definition \ref{Def-par-spec} from spaces to stacks.

\begin{definition}
[Smooth parametrized spectrum]
\label{def-2nd}
A \emph{smooth parametrized prespectrum} is a collection of morphisms 
$
p_n:E_n\to M
$
between smooth stacks in $\sh_{\infty}(\mathscr{M}{\rm f})$, $n\in \ZZ$, 
with a choice of section, equipped with morphisms $\Sigma_ME_n\to E_{n+1}$, 
 making similar diagrams
as in Definition \ref{Def-par-spec} 
%
commute up to a choice of equivalence in $\sh_{\infty}(\mathscr{M}{\rm f})$. 
A smooth parametrized prespectrum $\{p_n:E_n\to M\}$ for which the adjoint maps 
$E_n\to \Omega_BE_{n+1}$ are equivalences is called  \emph{smooth parametrized spectrum}. 
\end{definition}

\begin{remark}
[Identifying the proper category as a setting]
{\bf (i)} Note that Definition \ref{def-2nd} is almost verbatim the same as Definition 
\ref{Def-par-spec}, the only difference being where the objects $E_n$ and $M$ live (i.e. smooth 
stacks instead of spaces).  In this context we still have a mapping spectrum \footnote{Note that the mapping 
spectra are not smooth or parametrized; they are ordinary topological spectra.} between two smooth spectra. 
The resulting structure is again a stable infinity category and we denote this category by 
$T(\sh_{\infty}(\mathscr{M}{\rm f}))$. 

\item  
{\bf (ii)} We will be most concerned with the case when $M$ is a smooth manifold. By the Yoneda embedding,
 every smooth manifold embeds as an object in $\sh_{\infty}(\mathscr{M}{\rm f})$ via its sheaf of smooth
  plots, i.e., the sheaf sending $M$ to the set of smooth maps $N\to M$, with $N$ any other manifold. 

\item 
{\bf (iii)} One might wonder why the seemingly complicated $\infty$-category $T(\sh_{\infty}(\mathscr{M}{\rm f}))$ 
is necessary to work in. In particular, one might think that working with the more familiar category of sheaves of chain 
complexes should be more transparent. Note, however, that we are naturally led to the $\infty$-\emph{topos} 
$T(\sh_{\infty}(\mathscr{M}{\rm f}))$ for two reasons. First, passing to sheaves is necessary to capture the 
geometry of the de Rham complex (which is a crucial ingredient in defining Deligne cohomology). Second, 
in contrast to the category of sheaves of chain complexes, the axioms of an $\infty$-topos (in particular descent) 
make it a convenient setting to talk about bundles.
\end{remark}

We have the following natural definition for a smooth, locally trivial bundle of spectra.

\begin{definition}
[Smooth bundle of spectra]
Let $M$ be a smooth manifold. A \emph{smooth bundle of spectra $\pi:E\to M$ over $M$ with fiber the sheaf of 
spectra $\mathscr{R}$} is an object in $T(\sh_{\infty}(\mathscr{M}{\rm f}))_M$ satisfying the same properties 
as in Definition \ref{Def-bun-spec} with $M$ replacing $X$.
%
%
%
\end{definition}
We now wish to focus our scope to the case of periodic Deligne cohomology. Consider 
the sheaf of ring spectra given by applying the Eilenberg-MacLane functor $\H$ to the 
sheaf of chain complexes ${\mathcal D}(\op{ev})$ and ${\mathcal D}(\op{odd})$. In Section \ref{Sec-ccoh}, 
we saw that this ring spectrum represents periodic Deligne cohomology, in the sense that 
$$
\widehat{H}^{\op{ev}}(M;\ZZ[u,u^{-1}])\cong \pi_0\map\big(M;\H({\mathcal D}(\op{ev}))\big)
\quad 
\text{and}
\quad
\widehat{H}^{\op{odd}}(M;\ZZ[u,u^{-1}])\cong \pi_0\map\big(M;\H({\mathcal D}(\op{odd}))\big)\;.
$$
We would like to identify a large class of twists for this theory. To this end, let us consider the stack of twists 
in diagram \eqref{Tw-diag} with $\widehat{\mathscr{R}}$ the periodic differential ring spectrum given 
by both $\H({\mathcal D}(\op{ev}))$ and $ \H({\mathcal D}(\op{odd}))$, separately. At first, it might appear that 
we would get two stacks of twists corresponding to both the even and odd degrees; however, this is not the case.
\begin{proposition}
[Equivalence of stacks of even and odd twists for periodic Deligne cohomology]
\label{D-plus-D-minus}
We have a canonical equivalence of smooth stacks
$$
\widehat{{\rm Tw}}_{\H({\mathcal D}(\op{odd}))}\simeq \widehat{{\rm Tw}}_{\H({\mathcal D}(\op{ev}))}\;,
$$
induced by shifting both the ring spectrum $\H(\ZZ[u,u^{-1}])$ 
and the invertible periodic de Rham complex $\Omega^*[u,u^{-1}]$ up by one degree each.
\end{proposition}
\theproof
It is clear formally that shifting a module spectrum $\mathscr{R}_{\tau}$ up by 
one degree is a module spectrum over the ring spectrum $\mathscr{R}$, i.e., 
the module maps $ \mu: \mathscr{R}^m \wedge \mathscr{R}^n_h \to \mathscr{R}^{m+n}_h$ give rise to maps 
$\mathscr{R}^m \wedge \mathscr{R}^{n+1}_h \to \mathscr{R}^{m+n +1}_h$.
 Similarly, shifting a K-flat invertible module $\mathcal{L}$ is again a K-flat invertible module. 
 Moreover, given any equivalence
$$
\H(\mathcal{L})\simeq \mathscr{R}_{\tau}\wedge \H\RR\;,
$$
we get a corresponding equivalence at the level of the shifts. By the universal property 
of the pullback, we have an induced map at the level of the twists. For smooth periodic 
Deligne cohomology this takes the form
$$
\widehat{{\rm Tw}}_{\H({\mathcal D}(\op{odd}))}\longrightarrow \;
 \widehat{{\rm Tw}}_{\H({\mathcal D}(\op{ev}))}\;.
$$
It is immediate that this map admits an inverse induced by shifting down. 
\endofproof

Proposition \ref{D-plus-D-minus} implies that we do not have to consider the even and 
odd degrees separately, but we can view a given twist as corresponding to either spectrum. 
Henceforth, we will only refer to \emph{the} stack of twists of periodic Deligne cohomology 
and denote the stack simply by $\widehat{\rm Tw}$.

\begin{remark}
[Chern-Simons type hierarchy of trivilizations of {\v C}ech-Deligne cocycles]
\label{rem-CS}
Consider a {\v C}ech-Deligne cocycle $\eta=(\eta^{(0)},\eta^{(1)},
\eta^{(2)},\cdots,\eta^{(2k+1)})$ 
on a smooth manifold $M$, where $\eta^{(i)}$ is the cocycle data on the $i$-fold intersection, i.e.,
 $\eta^{(0)}$ is a $2k$-form defined on open sets, $\eta^{(1)}$ is a $(2k-1)$ form defined on intersections,  
 etc.  To this we associate automorphisms and higher automorphisms of the periodic de Rham complex 
 $\Omega^*[u,u^{-1}]$ on a smooth manifold $M$. More precisely, 
we associate to such a cocycle the automorphisms
\begin{align*}
 \ch^{(0)}(\eta)\wedge (-)&=e^{\eta^{(0)}} \wedge (-)=1+\eta^{(0)}\wedge (-)+
 \tfrac{1}{2!}\eta^{(0)}\wedge \eta^{(0)}\wedge (-) +
\tfrac{1}{3!}\eta^{(0)}\wedge \eta^{(0)}\wedge \eta^{(0)}\wedge (-)+\cdots \;,
\\
{\rm CS}^{(1)}(\eta)\wedge (-) &=\eta^{(1)}+
\tfrac{1}{2!}\eta^{(1)}\wedge d\eta^{(1)}\wedge(-)+
\tfrac{1}{3!}\eta^{(1)}\wedge d\eta^{(1)}\wedge 
d\eta^{(1)}\wedge(-)+\cdots \;,
\\
\vdots
\end{align*}
where we have a primary invariant  $\ch^{(0)}(\eta)$, 
then a secondary invariant ${\rm CS}^{(1)}(\eta)$ for the latter, 
then a tertiary invariant for the latter, and a similar pattern in higher degrees, 
obtained by using the various cocycle data for $\eta$. 
This assignment is 
explained in the proof of \cite[Theorem 17]{GS5} and in the discussion leading up to that 
theorem. In particular, taking $k=1$, we get $\eta=(B_\a, A_{\a \b}, f_{\a \b \gamma},
n_{\a \b \gamma, \delta})$, which is the {\v C}ech-Deligne cocycle 
corresponding to the `standard' gerbe with connection (see Example \ref{bry gerbe}), 
encoding  a twist.  These represent the homotopies, homotopies between homotopies, etc., respectively, 
in diagram \eqref{cdiag-sim} of Example \ref{Ex-shCS} below. 
The  expression for ${\rm CS}^{(1)}(\eta)$ is a sum of 
higher product abelian Chern-Simons theories, in the sense of \cite{FSS1}. 
The next terms (not explicitly recorded for brevity) correspond to tertiary and higher 
structures, in the sense of \cite{FSS1}\cite{tert}.
\end{remark}

\begin{proposition}
[Twisting periodic Deligne cohomology by odd degree gerbes with connection]
\label{Prop-Del-gerbe}
Let $\widehat{{\rm Tw}}(M)$ denote the stack of twists, evaluated on a smooth manifold $M$. Then 
every {\v C}ech-Deligne cocycle of degree $2k+1$ defines a twist of periodic Deligne cohomology. 
In fact, there is a morphism of smooth stacks
$$
\BB^{2k}U(1)_{\nabla}\longrightarrow \widehat{{\rm Tw}}\;,
$$
refining the map
$
K(\ZZ,2k+1)\to B{\rm GL}_1(\H\ZZ[u,u^{-1}])
\into {\rm Pic}^{\rm top}_{\H\ZZ[u,u^{-1}]}
$.
\end{proposition}\label{twperd diff}
\theproof
The stack $\BB^{2k}U(1)_{\nabla}$ fits into the Cartesian square \eqref{stack gpull}. 
This pullback in smooth stacks can be computed by the stackification of the corresponding pullback 
in prestacks,  which is computed objectwise. Moreover,  a morphism of prestacks into a stack is, 
equivalently, a morphism of stacks out of the stackification.\footnote{This follows from the 
adjunction $i\vdash L$, with $L$ the stackification functor and $i$ the inclusion functor.}
  It, therefore, suffices to construct the map objectwise out of the three stacks $ \Omega^{2k+1}_{\rm cl}$, 
  $\BB^{2k+1}\ZZ$ and $\Omega_{\rm cl}^{\leq 2k+1}$ and for every (1-)homotopy filling the diagram, 
  a corresponding homotopy filling diagram \eqref{Tw-diag}. 

To that end, fix an arbitrary manifold $M$ and define the three maps 
\begin{align*}
&\;\;\;\xymatrix{\Omega_{\rm cl}^{2k+1}(M) \ar[r] &  {\rm Pic}^{\rm form}(M)}
\\
&\xymatrix{\BB^{2k+1}\RR(M)  \ar[r]& \BB {\rm GL}_1(\H\RR[u,u^{-1}])(M)
\;\ar@{^{(}->}[r]& {\rm Pic}^{\rm dR}(M)}
\\
&\xymatrix{\BB^{2k+1}\ZZ(M) \ar[r]& \BB {\rm GL}_1(\H\ZZ[u,u^{-1}])(M)
\; \ar@{^{(}->}[r]& {\rm Pic}^{\rm top}(M)}
\end{align*}
as follows. The first map sends a closed odd-degree form to the invertible module over the periodic de 
Rham complex, $(\Omega^*[u,u^{-1}](M),d_{H})$, where the differential $d_{H}=d+H\wedge $ 
acts on a differential form as
\(
\label{dhomega}
d_{H}(\omega)=d_{H}(\omega_0+\omega_2+\cdots)=d\omega_0+d\omega_2+\cdots + 
(H\wedge \omega_0+d\omega_{2k})+\cdots \;.
\)
The second map is induced by the canonical inclusion map $K(2k+1,\RR)\into 
B{\rm GL}_1(\H\RR[u,u^{-1}])$, while the third map is similarly induced by 
$K(2k+1,\ZZ)\into B{\rm GL}_1(\H\ZZ[u,u^{-1}])$.

A homotopy filling the diagram \eqref{stack gpull} can be identified via the Dold-Kan correspondence as an 
element $\eta$ in degree 1 of the total complex of the {\v C}ech double complex  
$
{\rm Tot}\big(\{U_{\alpha}\},\Omega^{\leq 2k+1}_{\rm cl}\big)
$,
with $\{U_{\alpha}\}$ a good open cover of $M$, subject to the following condition:
We require that $D(\eta)=H-h$, with $H$ the globally defined closed differential form of degree $2k+1$ 
determined in $\Omega^{2k+1}_{\rm cl}(M)$ and $h$ the real-valued {\v C}ech cocycle of degree $2k+1$ 
representing an element in $\BB^{2k+1}\ZZ(M) \into \BB^{2k+1}\RR(M)$ (see \cite{Cech} for more details 
on the latter stack). To the homotopy $\eta$, we need to construct a corresponding equivalence
$$
\H(\Omega^*[u,u^{-1}],d_{H})\simeq \H\RR[u,u^{-1}]_h\;,
$$
where $\H\RR[u,u^{-1}]_h$ is the locally constant sheaf of spectra corresponding to an invertible module over 
$\H\RR[u,u^{-1}]$.  The former is classified by the map 
$$
\xymatrix{
h:M\ar[r] & K(\ZZ,2k+1)\ar[r]&
 K(\RR,2k+1)\;\ar@{^{(}->}[r] & {\rm Pic}^{\rm dR}_{\H\RR[u,u^{-1}]}
}\;.
$$ 
As described in \cite{GS5}, the locally defined form $\ch^{(0)}(\eta)$ 
from Remark \ref{rem-CS}
defines a local trivialization 
$$
\ch^{(0)}(\eta):(\Omega^*[u,u^{-1}],d_{H})\vert_{U_{\alpha}}\overset{\simeq}{\longrightarrow} 
\Omega^*[u,u^{-1}]\vert_{U_{\alpha}}\;.
$$ 
The higher $CS^{(i)}$'s in Remark \ref{rem-CS} correspond to automorphisms on intersections and higher 
automorphisms on higher intersections. As in the proof of \cite[Theorem 17]{GS5}, this cocycle determines 
an edge in ${\rm Pic}^{\rm dR}$ connecting $\H(\Omega^*[u,u^{-1}],d_{H})$ and $\H\RR[u,u^{-1}]_h$. 
It is clear from the construction that this map refines the inclusion of $K(\ZZ,2k+1)$ into the topological twists. 
\endofproof

In Section \ref{Sec Per Z}, we described a universal bundle of spectra which classifies bundles of spectra over a space $X$ 
with a  $K(\ZZ,2k)$-structure prescribed by a twist $h:X\to K(\ZZ,2k+1)$. There is a similar universal bundle over 
the stack of twists \eqref{Tw-diag}, described in \cite{GS5}, which classifies twisted differential cohomology theories via pullback. 
Descent allows us to glue together the bundle via local trivializations. 
As a fundamental example, we consider the following.
\begin{example}
[Higher-twisted differential forms as sections of a smooth bundle of spectra]  
\label{Ex-shCS}
Consider the sheaf of complexes $(\Omega^*[u,u^{-1}],d_{H})$ on a smooth manifold $M$, which is degreewise identical 
to the periodic complex of forms, but which is equipped with the differential $d_{H}:=d+H\wedge$, acting by \eqref{dhomega}. 
Here, $H$ a closed form of degree $2k+1$. 
Applying the Eilenberg-MacLane functor ${\H}$ to $(\Omega^*[u,u^{-1}],d_{H})$ gives a sheaf of spectra on $M$. 
Now $(\Omega^*[u,u^{-1}],d_{H})$ is an invertible module over $(\Omega^*[u,u^{-1}], d)$ which is locally 
equivalent (by the Poincar\'e Lemma) to the constant sheaf $\RR[u,u^{-1}]$. Thus 
$\H(\Omega^*[u,u^{-1}],d_{H})$ gives a sheaf of spectra which is a module over $\H\RR[u,u^{-1}]$. 
Pulling back  the universal bundle of spectra
$$
\lambda\longrightarrow \; {\rm Pic}_{\H\RR[u,u^{-1}]}^{\rm dR}
$$ 
(see \cite{GS5} for this construction) by the map $\tau:M\to {\rm Pic}_{\H\RR[u,u^{-1}]}^{\rm dR}$, 
which picks out the twisted sheaf of spectra $\H(\Omega^*[u,u^{-1}],d_{H})$, gives a smooth bundle of spectra 
$E\to M$. The sheaf of local sections of the latter evaluated on $U$ is, by definition, 
$\H(\Omega^*[u,u^{-1}],d_{H})(U)$. Choose local potentials $B_{\alpha}$ for $H$ on each element of
 a good open cover $\{U_{\alpha}\}$ of $M$ (i.e. $dB_{\alpha}=H$). Then, on each patch $U_{\alpha}$, we have 
 quasi-isomorphisms of sheaves of complexes
$$
e^{B_{\alpha}}\wedge :(\Omega^*[u,u^{-1}],d_{H})\vert_{U_{\alpha}}\overset{\simeq}{\longrightarrow} \Omega^*[u,u^{-1}]\vert_{U_{\alpha}}
$$
which send a local section $\omega$ to the wedge product with the formal exponential
$$
e^{B_{\alpha}}=1+B_{\alpha}+\tfrac{1}{2!}B_{\alpha}^2+\cdots.
$$
These quasi-isomorphisms correspond to local trivializations
$$
e^{B_{\alpha}}\wedge :E\vert_{U_{\alpha}}\longrightarrow
 \H(\Omega^*[u,u^{-1}])\times U_{\alpha}\;.
$$
In fact, a choice of representative of $H$ in the {\v C}ech-de Rham double complex 
gives rise to a homotopy commutative diagram 
\(
\label{cdiag-sim}
\hspace{-.6cm}
\xymatrix@=2em{
\hdots \ar@<.15cm>[r] \ar@<.05cm>[r] \ar@<-.05cm>[r] \ar@<-.15cm>[r] &  \coprod_{\alpha\beta\gamma} 
\H(\Omega^*[u,u^{-1}])\times  U_{\alpha\beta\gamma}\ar@<.1cm>[r] \ar[r] \ar@<-.1cm>[r] & \coprod_{\alpha\beta} \H(\Omega^*[u,u^{-1}])\times U_{\alpha\beta}\ar@<.05cm>[r] \ar@<-.05cm>[r] 
& \coprod_{\alpha}\H(\Omega^*[u,u^{-1}])\times U_{\alpha}
}
\)
where the simplicial maps at each stage are determined by the {\v C}ech-de Rham data for $H$. For example, choosing 
a differential form $A_{\alpha\beta}$ on intersections, satisfying $dA_{\alpha\beta}=B_{\alpha}-B_{\beta}$, gives
 rise to a homotopy commutative diagram 
$$
\xymatrix@R=1.5em{
\H(\Omega^*[u,u^{-1}])\times  U_{\alpha\beta}\ar[rr]^-{{\rm id}} \ar[dr]_{e^{-B_{\alpha}}} 
& &\H( \Omega^*[u,u^{-1}])\times  U_{\alpha\beta}
\\
& E\vert_{U_{\alpha\beta}}\ar[ru]_-{e^{B_{\beta}}}&
}\;,
$$
with homotopy given by wedge product with the abelian Chern-Simons form (cf. Remark 
\ref{rem-CS})
$$
CS(A_{\alpha\beta})=A_{\alpha\beta}+
\tfrac{1}{2!}A_{\alpha\beta}\wedge dA_{\alpha\beta}
+\tfrac{1}{3!}A_{\alpha\beta}\wedge dA_{\alpha\beta}\wedge dA_{\alpha\beta}
+\cdots
$$ 
so that
$$
dCS(A_{\alpha\beta})=e^{dA_{\alpha\beta}}=e^{B_{\alpha}-B_{\beta}}=e^{B_{\alpha}}e^{-B_{\beta}}\;.
$$
By descent, one concludes that $E$ is, in fact, a colimit over this diagram. Notice also 
that, since the global sections of $E$ are $\H(\Omega^*[u,u^{-1}](M),d_{H})$, we 
immediately have that the twisted cohomology represented by the bundle $H$ is 
the twisted de Rham cohomology of $M$.
\end{example}

\begin{remark}
[The twisted periodic Deligne complex as sections of a smooth bundle of spectra] 
Example \ref{Ex-shCS} shows that, secretly, elements in the twisted de Rham complex are sections 
of a bundle of spectra over $M$. That example can be adapted to the periodic Deligne 
complex by noting that if the de Rham twist $H$ has integral periods, then the corresponding topological
 twist $h:M\to K(\RR,2k+1)$ factors (up to homotopy) through $K(\ZZ,2k+1)$. In this case, the triple 
$\big(H\ZZ[u,u^{-1}]_h,t,(\Omega^*[u,u^{-1}],d_{H})\big)$, with 
$$
t:\H\RR[u,u^{-1}]_h\simeq \H(\Omega^*[u,u^{-1}],d_{H})
$$
the twisted de Rham equivalence, gives a twist of periodic Deligne cohomology on $M$, hence a map
$\tau:M\to \widehat{\rm Tw}$. Pulling back by the universal bundle $\hat{\lambda}\to \widehat{\rm Tw}$
 gives similar data to the ones in Example  \ref{Ex-shCS}.

\end{remark}

\begin{definition}
[Twisted periodic Deligne cohomology]
Let $H$ be a closed differential form of odd degree which has integral periods. Then 
$H$ can be lifted to a map 
$$
\hat{h}:M\longrightarrow \BB^{2k}U(1)_{\nabla}\;.
$$
According to Proposition \ref{D-plus-D-minus}, we can regard this as either twisting $\H({\mathcal D}(\op{ev}))$ or
$\H({\mathcal D}(\op{odd}))$. Let ${\mathcal Z}_{\hat{h}}^{\op{ev}}\to M$ and ${\mathcal Z}_{\hat{h}}^{\op{odd}}\to M$ 
be the corresponding smooth bundles of spectra. We define the \emph{twisted periodic Deligne cohomology} to be the 
homotopy classes of sections of the corresponding bundle ${\mathcal Z}_{\hat{h}}^{\op{ev}/\op{odd}}\to M$, i.e.,
$$
\widehat{H}^{\op{ev}}(M;\hat{h}):=\pi_0\Gamma(M;{\mathcal Z}_{\hat{h}}^{\op{ev}})
\qquad
\text{and} 
\qquad
\widehat{H}^{\op{odd}}(M;\hat{h}):=\pi_0\Gamma(M;{\mathcal Z}_{\hat{h}}^{\op{odd}})\;.
$$
\end{definition}

\medskip
Just as one can define tensor product and direct sum of vector bundles, one can similarly define the wedge product and 
smash product of bundles of spectra (see \cite{GS5} for the definition of the smash product; the wedge product is defined 
similarly). In the present case, we have a $\ZZ/2$-graded bundle of spectra $E^{\op{ev}}\vee E^{\op{odd}}\to M$. The 
local sections of this bundle are given by evaluating the wedge product of spectra
$\H({\mathcal D}(\op{ev}))\vee \H({\mathcal D}(\op{odd}))\simeq \H({\mathcal D}(\op{ev})\oplus {\mathcal D}(\op{odd}))$ 
on open subsets $U\subset M$. Given the multiplicative structure of periodic Deligne cohomology
from Section \ref{Sec-ring}, we have the following.

\begin{proposition}
[Module structure of twisted periodic Deligne cohomology]
The sheaf of sections of the wedge bundle $E^{\op{ev}}\vee E^{\op{odd}}\to M$ is a module spectrum over the 
sheaf of ring spectra given by the wedge 
$\H({\mathcal D}(\op{ev}))\vee \H({\mathcal D}(\op{odd}))$. The module action descends to a map
$$
\xymatrix{
\mu:\widehat{H}^{\op{ev}/\op{odd}}(M;\ZZ[u,u^{-1}])\otimes \widehat{H}^{\op{ev}/\op{odd}}(M;\hat{h})
\ar[r] &
\widehat{H}^{\op{ev}/\op{odd}}(M;\hat{h})}\;,
$$
turning $\widehat{H}^{\op{ev}/\op{odd}}(M;\hat{h})$ into a module over the 
superalgebra $\widehat{H}^{\op{ev}/\op{odd}}(M;\ZZ[u,u^{-1}])$.
\end{proposition}

\medskip
\subsection{Properties of twisted periodic smooth Deligne cohomology}
\label{Sec-cech} 

In this section, we list some of the properties of twisted periodic Deligne cohomology. Most of these properties are 
similar to those discussed in Proposition \ref{prop-tpic}. 
A main point of contrast to emphasize here is that this theory will not satisfy
the homotopy invariance axiom, as is the case for any differential cohomology theory (see \cite{BS}). 

\begin{proposition}
[Properties of twisted periodic Deligne cohomology]
\label{pr twdel}
Let $M$ be a smooth manifold and fix a twist $\hat{h}:M\to \BB^{2k}U(1)_{\nabla}$. Consider the category of such 
pairs $(M,\hat{h})$, with morphisms $f:(M,\hat{h})\to (M, \hat{\ell})$ given by smooth maps $f:M\to N$ such that 
$[f^*\hat{\ell}]=[\hat{h}]$. The assignment $(M,\hat{h})\mapsto \widehat{H}^*(M;\hat{h})$ satisfies the following 
properties:
\begin{enumerate}[{\bf (i)}]
\item $\widehat{H}^*(M;\hat{h})$ is functorial with respect to the maps $f:(M,\hat{h})\to (N,\hat{\ell})$.

\item The functor $\widehat{H}^*(-;\hat{h})$ satisfies the Eilenberg-Steenrod axioms (modulo the dimension 
axiom and homotopy invariance!) for a reduced cohomology theory. In particular, we have a Mayer-Vietoris 
sequence which takes the form
\begin{equation*}
\xymatrix{
  &
 ... \ar[r] &  
H_{\RR/\ZZ}^{*-2}(U;h)
\oplus
H_{\RR/\ZZ}^{*-2}(V;h)
 \ar[r]
 & 
H_{\RR/\ZZ}^{*-2}(U\cap V;h)
 \ar`r[d]`[l]`[llld]`[dll][dll]  \\
  &\widehat{H}^{*}(M;\hat{h}) \ar[r] & 
\widehat{H}^{*}(U;\hat{h})\oplus\widehat{H}^{*+1}(V;\hat{h})\ar[r]& 
\widehat{H}^{*}(U\cap V ;\hat{h}) \ar`r[d]`[l]`[llld]`[dll][dll] 
\\ 
 &H^{*+1}(M;h) \ar[r] &  H^{*+1}(U;h)
 \oplus H^{*+1}(V;h) \ar[r]& ...
 }
\end{equation*}

\item For $\hat{h}:M\to \BB^{2k}U(1)_{\nabla}$ a trivial twist (i.e. $\hat{h}\simeq \ast$ in smooth stacks) 
we have an isomorphism
\(\label{iso twuntw}
\widehat{H}^*(M;\hat{h})\cong \widehat{H}^*(M;\ZZ[u,u^{-1}])\;.
\)
Even more strongly, we still have an isomorphism \eqref{iso twuntw} if just the underlying 
topological twist $h:M\to K(\ZZ,2k+1)$ is trivial.

\end{enumerate}
\end{proposition}

\theproof
{\bf (i)}   Given a map $f:M\to N$ satisfying the desired compatibility, we have an 
induced double pullback diagram 
\(
\label{higher bundle czh}
\xymatrix@R=1.5em{
 {\mathcal Z}_{f^*\hat{h}}\ar[r]\ar[d] & {\mathcal Z}_{\hat{h}} \ar[rr]\ar[d] &&  \hat{\lambda}\ar[d]
\\
M \ar[r]^f & N\ar[r]^-{\hat{h}} & \BB^{2k}U(1)_{\nabla} \ar[r] & \widehat{\rm Tw}
\;,
}
\)
which gives the identification ${\mathcal Z}_{f^*\hat{h}}\simeq {\mathcal Z}_{\hat{\ell}}$.
As a consequence, we have an induced morphism of sections $f^*:\Gamma(N;{\mathcal Z}_{\hat{h}})\to \Gamma(M;{\mathcal Z}_{\hat{\ell}})$. Passing to homotopy groups yields a map $f^*:\widehat{H}^*(N;\hat{\ell})\to \widehat{H}^*(M;\hat{h})$. 

\medskip
\noindent {\bf (ii)} We now verify the applicable Eilenberg-Steenrod axioms. 

\vspace{0.2cm}

\noindent {\it Additivity}. 
Let $M=\coprod_{\alpha}M_{\alpha}$ with each $M_{\alpha}$ a smooth manifold. A map 
$\hat{h}:M\to \BB^{2k}U(1)_{\nabla}$ is equivalently a collection of maps 
$h_{\alpha}:M_{\alpha}\to \BB^{2k}U(1)_{\nabla}$. Then the spectrum of sections 
$\Gamma(M,{\mathcal Z}_{\hat{h}})$ splits as a
product $\prod_{\alpha}\Gamma(M_{\alpha},{\mathcal Z}_{\hat{h}_{\alpha}})$. 
Since taking homotopy groups commutes with products,  we have an isomorphism
$$
\widehat{H}^*(M,\hat{h})\cong \prod_{\alpha}\widehat{H}^*(M_{\alpha},\hat{h}_{\alpha})\;.
$$
\noindent {\it Exactness}. This follows verbatim as in the proof of Proposition \ref{prop-tpic}, with 
the space $X$ replaced by a smooth manifold $M$, $A\subset M$ a submanifold, and 
the map $h:X\to K(\ZZ,2k+1)$ replaced by the refinement $\hat{h}:M\to \BB^{2k}U(1)_{\nabla}$. 


\vspace{2mm}
\noindent {\bf (iii)} Finally, if the twist $\hat{h}:M\to \BB^{2k}U(1)_{\nabla}$ is topologically trivial, i.e., 
its geometric realization $h:\vert M\vert \simeq M\to \vert \BB^{2k}U(1)_{\nabla}\vert\simeq K(\ZZ,2k+1)$ 
is homotopic to a the constant map induced by $0\to \ZZ$. In this case, the underlying twisted spectrum 
$\H\ZZ[u,u^{-1}]_{h}$ is equivalent to $\H\ZZ[u,u^{-1}]$ and we have the diagram 
$$
\xymatrix@R=1.8em{
\H(\ZZ[u,u^{-1}] )\ar[rr]^-\simeq \ar[d]_-{\wedge \H \RR} &&
\H(\ZZ[u,u^{-1}])_h \ar[d]^-{\wedge \H \RR}
\\
\H(\RR[u,u^{-1}]) \ar[d]_{\simeq}^{c} \ar[rr]^-\simeq &&
\H(\RR[u,u^{-1}])_h \ar[d]^{\simeq}_{t} 
\\
\H(\Omega^*[u,u^{-1}]) \ar[rr]^-{\simeq } && \H(\Omega^*[u,u^{-1}] ,d_{H}) \;,
}
$$
where the bottom equivalence depends on a choice of homotopy inverse for $c$ and is defined as the obvious 
composition in the diagram. By the basic properties of the functor $\H$ (see \cite[pp. 17-18]{BN} for discussion), 
the existence of the bottom equivalence implies that $\Omega^*[u,u^{-1}]$ and $(\Omega^*[u,u^{-1}] ,d_{H})$
 are connected by a zig-zag of quasi-isomorphisms. This is manifestly the data  needed to define an equivalence 
 in the $\infty$-groupoid $\widehat{\rm Tw}(M)$.
\endofproof

\section{The spectral sequence $\widehat{\rm AHSS}_h$ 
and examples}
\label{Sec-tAHSS}

In this section, we apply the twisted Atiyah-Hirzebruch spectral sequence 
(both the classical \cite{Ro1}\cite{Ro2}
\cite{AS} and the differential refinement 
 \cite{GS5}) to calculate the twisted periodic integral and Deligne cohomology of spheres.

\subsection{The spectral sequence in twisted periodic smooth Deligne cohomology}
\label{tAHSS}

In \cite{AS}, the first nonvanishing differential for the twisted AHSS (applied to $K$-theory) 
on a space $X$ was identified by observing that the only degree three increasing operations for 
spaces equipped with maps $X\to K(\ZZ,3)$ are given by the cohomology group
$$
H^{n+3}(K(\ZZ,n)\times K(\ZZ,3))\cong H^{n+3}(K(\ZZ,n))\oplus H^{n+3}(K(\ZZ,3))\oplus \ZZ\;.
$$
The third factor on the right hand side is generated by the product of the generators 
for $H^{n}(K(\ZZ,n))$ and $H^{3}(K(\ZZ,3))$. From this, one deduces that 
$$
d_3(x)=Sq^3_{\ZZ}(x)- [h]\cup x\;,
$$
with $[h]$ the twisting integral class and $Sq^3_{\ZZ}$ the third integral Steenrod square, 
which comes from the untwisted AHSS for $K$-theory \cite{AH}. 

\medskip
The situation for periodic integral cohomology is much easier since the untwisted differentials vanish and 
since one is able to compare easily with rational cohomology. Considering again a degree three twist $h$, 
we then find that \footnote{This can be deduced, for example, from the $K$-theory differential and the 
fact that (on spheres) the Chern character lands in integral cohomology.}
$$
d_3(x)=- [h] \cup x\;.
$$
The same argument applies not only in the degree three case, but also in higher odd degrees.
This is due to the fact that for spaces equipped with maps $X\to K(2k+1;\ZZ)$, we again have the identification
$$
H^{n+2k+1}(K(\ZZ,n)\times K(\ZZ,2k+1))\cong H^{n+2k+1}(K(\ZZ,n))\oplus H^{n+2k+1}(K(\ZZ,2k+1))\oplus \ZZ\;,
$$
with the last factor being generated by the product of the generator of $H^{n}(K(\ZZ,n))$ and 
the generator of $H^{2k+1}(K(\ZZ,2k+1))$. We, therefore, have the following.
\begin{proposition}
[First differential for $\op{AHSS}_h$ for twisted periodic integral cohomology]
Let $h:X\to K(\ZZ,2k+1)$ be a twist of periodic integral cohomology. Then the first 
nonvanishing differential in the associated AHSS occurs on the $E_{2k+1}$-page and is given by
$$
d_{2k+1}(x)=- [h] \cup x\;.
$$
\end{proposition}

We will illustrate this in Examples \ref{Ex-Z-Sev} and \ref{Ex-Z-Sodd} below. 

\medskip
In \cite{GS5}, we developed an AHSS for twisted \emph{differential} cohomology theories, in turn generalizing 
that of a differential theory \cite{GS3}. In the case of periodic Deligne cohomology, there are two spectral sequences 
corresponding to the even and odd degrees (separately). 

\begin{lemma}
[The $E_2$-page for even degrees in $\widehat{\rm AHSS}_{\hat h}$ for twisted periodic Deligne cohomology]
The $E_2$-page for the even case looks as follows
{\small 
\(
\label{E2-Kth}
\begin{tikzpicture}
  \matrix (m) [matrix of math nodes,
    nodes in empty cells,nodes={minimum width=5ex,
    minimum height=4ex,outer sep=-3pt}, 
    column sep=1ex,row sep=1ex]{
                &      &     &     & 
                \\         
            1 &     &   &   &
            \\
             0     &   \Omega^{{\rm ev}}_{{d_H}\mbox{-}{\rm cl}, \ZZ}(M) &  &   &
             \\           
               -1\  \   &  & H^1(M;U(1)) &  H^2(M;U(1))     & 
                \\
              -2\ \  &   &    &  0 & 0 & 0
               \\
               -3\ \ &  & & &  & H^4(M;U(1))
               \\
               -4\ \ &   &    &   &  &  & 
               \\
                &      &     &     & \\};
                
  \draw[-stealth] (m-3-2.east) -- node[above]{\footnotesize $d_2$} (m-4-4.north west);
   \draw[-stealth] (m-4-3.south east) -- (m-5-5.west);
      \draw[-stealth] (m-5-4.south east) -- (m-6-6.north west);
        \draw[-stealth] (m-4-4.south east) -- (m-5-6.north west);
\draw[thick] (m-1-1.east) -- (m-8-1.east) ;
\end{tikzpicture}
\)
}
where  
$\Omega^{{\rm ev}}_{{d_H}\mbox{-}{\rm cl}, \ZZ}(M)$
is the subgroup of those even forms on $M$ which are twisted-closed
and whose  degree zero component is given by an integer, i.e. 
$\omega=n_0+\omega_2+\omega_4+\cdots$. 
\end{lemma}

\begin{lemma}
[The $E_2$-page for odd degrees in $\widehat{\rm AHSS}_{\hat h}$ for twisted periodic Deligne cohomology]
The spectral sequence for the odd degrees looks as follows
{\small 
\(
\label{E2-Kth}
\begin{tikzpicture}
  \matrix (m) [matrix of math nodes,
    nodes in empty cells,nodes={minimum width=5ex,
    minimum height=4ex,outer sep=-2pt}, 
    column sep=1ex,row sep=1ex]{
                &      &     &     & 
                \\         
            1 &     &   &   &
            \\
             0     &   U(1)\times \Omega^{{\rm odd}}_{{d_H}\mbox{-}{\rm cl}}(M)&  &   &
             \\           
               -1\  \   &  & 0 &  0    & 
                \\
              -2\ \  &   &    &  0 & H^2(M;U(1)) & H^3(M;U(1)) 
               \\
               -3\ \ &  & & &  & 0
               \\
               -4\ \ &   &    &   &  &  H^4(M;U(1))& 
               \\
                &      &     &     & \\};
                
  \draw[-stealth] (m-3-2.south) -- ($(m-4-4.north west)!.4!(m-4-4.south west)$);
   \draw[-stealth] (m-4-3.south east) -- ($(m-5-5.north east)!0.12!(m-5-2.north west)$);
      \draw[-stealth] (m-5-4.south east) -- ($(m-6-6.north west)!.2!(m-6-6.south west)$);
        \draw[-stealth] (m-4-4.south east) -- (m-5-6.north west);
\draw[thick] (m-1-1.east) -- (m-8-1.east) ;
\end{tikzpicture}
\)
}
where $\Omega^{{\rm odd}}_{{d_H}\mbox{-}{\rm cl}}(M)$ is the group of twisted-closed odd forms
on $M$.
\end{lemma}
For twisted differential $K$-theory, in \cite{GS5} we identified the first nonzero differential 
in the spectral sequence as
$$
d_3(x)=\widehat{Sq}^3_{\ZZ}(x)+ [\hat{h}]\cup_{\rm DB} x\;,
$$
where $\widehat{Sq}^3_{\ZZ}$ is a torsion operation in differential cohomology inherited from 
$Sq^3$ (see \cite{GS2}), and $[\hat{h}]\cup _{\rm DB}(-)$ is the Deligne-Beilinson cup product 
operation. The same argument used in \cite[Proposition 25]{GS5}  applies to the case of differential 
refinements of the higher degree twists for periodic integral cohomology. As a result we have the following. 

\begin{proposition}
[First differential for $\widehat{\rm AHSS}_{\hat h}$ for twisted Deligne cohomology]
Let $h:M\to \BB^{2k}U(1)_{\nabla}$ be a twist of periodic Deligne cohomology.
 Then the differential in the associated AHSS on the $E_{2k+1}$-page
 \footnote{Note that there is also a differential on the $E_{2k}$-page; but we do not use this.}
  is given by
$$
d_{2k+1}(x)=- [\hat{h}]\cup_{\rm DB} x\;.
$$
\end{proposition}

We will illustrate this in Examples \ref{Ex-Del-Sev} and \ref{Ex-Del-Sodd} below.

\subsection{Examples via the spectral sequence} 
\label{Sec-Ex}

We now proceed with our examples illustrating the AHSS
that we developed in Section \ref{tAHSS} to both 
twisted periodic integral cohomology (Section \ref{Sec Per Z})
and  twisted periodic Deligne cohomology (Section \ref{Sec Higher}). 

\begin{example}[Twisted periodic integral cohomology of even spheres]
\label{Ex-Z-Sev}
For even spheres, the class of the twist vanishes for parity reasons. Therefore,
 the AHSS degenerates at the $E_{2}$-page and we immediately identify
$$
H^{\rm ev}(S^{2k};\ZZ[u,u^{-1}])=\ZZ\oplus \ZZ
\qquad \text{and} \qquad
H^{\rm odd}(S^{2k},\ZZ[u,u^{-1}])=0\;.
$$
\end{example}

\begin{example}[Twisted periodic integral cohomology of odd spheres]
\label{Ex-Z-Sodd}
For an odd-dimensional sphere $S^{2k+1}$, the only interesting twist occurs in degree $2k+1$. 
Consequently, the  only nonzero differential in the AHSS occurs on the $E_{2k+1}$-page, 
giving the sequence
$$
\xymatrix@=1.5em{
\ZZ\ar[rr]^-{-[h]\cup} && H^{2k+1}(S^{2k+1}:\ZZ)\cong \ZZ\ar[r] & 0
}\;.
$$
Thus, with ${h}$ also denoting the integer corresponding to the topological 
twist $h$, we get
$$
H^{\rm odd}(S^{2k+1};\ZZ[u,u^{-1}])\cong \ZZ/{h} 
\qquad 
\text{and}
\qquad
H^{\rm ev}(S^{2k+1};\ZZ[u,u^{-1}])\cong 0\;.
$$
\end{example}

The above examples illustrate the utility of the AHSS in computing 
the twisted integral cohomology. 
We will now extend these examples to the differential case. 
It turns out that for even spheres, we will find that the use of the Mayer-Vietrois sequence
is straightforward enough and efficient in this case.

\begin{example}[Twisted periodic Deligne cohomology of even spheres]
\label{Ex-Del-Sev}
Let $\hat{h}:M\to \BB^{2k}U(1)_{\nabla}$ be a twist for periodic Deligne cohomology. 
For parity reasons, the class of the underlying topological twist $h\in H^{2k+1}(S^{2k};\ZZ)$ 
vanishes. By property (iii) of Proposition \ref{pr twdel}, 
it follows that we have an isomorphism
$$
\widehat{H}^{\rm ev/odd}(S^{2k};\hat{h})\cong  \widehat{H}^{\rm ev/odd}(S^{2k};\ZZ[u,u^{-1}])
$$
with the underlying untwisted theory. We computed the corresponding groups  earlier 
in Example \ref{Ex-PD-S2k} (Section \ref{Sec-ring}), which immediately yields
 $$
\widehat{H}^{\rm ev}(S^{2k}; \hat{h})= 
\Omega^{\op{odd}}(S^{2k})/{\rm im}(d)\oplus \ZZ\oplus \ZZ
\qquad 
\text{and}
\qquad
\widehat{H}^{\rm odd}(S^{2k},\hat{h})= 
\Omega^{\op{ev}}(S^{2k})/\Omega^{\op{ev}}_{\rm cl,\ZZ}(S^{2k})\;.
$$
\end{example}

\medskip 
The case of odd spheres is more involved.

\begin{example}[Twisted periodic Deligne cohomology of odd spheres]
\label{Ex-Del-Sodd}
For the odd spheres, the only interesting twist are the differential refinements of the topological 
twists $[h]\in H^{2k+1}(S^{2k+1};\ZZ)$. Choose such a differential refinement 
$\hat{h}:M\to \BB^{2k}U(1)_{\nabla}$. Then the spectral sequence has one nontrivial differential
$$
\xymatrix{
d_{2k+1}:U(1)\times \Omega^{\rm odd}_{d_H\mbox{-}{\rm cl}}(S^{2k+1})
\ar[r] & \;  U(1)\cong H^{2k+1}(S^{2k+1};U(1))}\;,
$$
occurring on the $(2k+1)$-page. Here $\Omega^{\rm odd}_{d_H\mbox{-}{\rm cl}}$ 
denotes those odd forms which are closed under the twisted differential $d_H$. 
The restriction to the factor $U(1)$ is given by the Deligne-Beilinson cup product with $[\hat{h}]$. 
This can be computed as follows. As above, let ${h}$ be the integer representing the topological class 
$h\in \ZZ$. For $\theta\in U(1)$, we have $[\hat{h}]\cup_{DB}\theta=h\theta$. The kernel of 
$d_{2k+1}$ restricted to this factor is the subgroup of $h$-roots of unity which is isomorphic to 
$\ZZ/h$. Since the map $\theta\mapsto h\theta$ is surjective 
the First Isomorphism Theorem implies that the factor 
$\Omega^{\rm odd}_{d_H\mbox{-}{\rm cl}}(S^{2k+1})$ is killed by $d_{2k+1}$. 

It remains to solve the extension problem 
$$
\xymatrix{
0\ar[r] & \ZZ/h\ar[r] & \widehat{H}^{\rm odd}(S^{2k+1};\hat{h})
\ar[r] & \Omega^{\rm odd}_{d_H\mbox{-}{\rm cl}}(S^{2k+1})
\ar[r] &0}\;.
$$
Now for any abelian group $A$ and any divisible group $B$ 
 the Ext group $\op{Ext}(A, B)$ vanishes. 
Since the group $\Omega^{\rm odd}_{d_H\mbox{-}{\rm cl}}(S^{2k+1})$ is divisible, 
we then have
$$
{\rm Ext}^1\big(\ZZ/h,\Omega^{\rm odd}_{d_H\mbox{-}{\rm cl}}(S^{2k+1})\big)
\cong
 \Omega^{\rm odd}_{d_H-{\rm cl}}(S^{2k+1})/h\Omega^{\rm odd}_{d_H\mbox{-}{\rm cl}}(S^{2k+1}) 
 \cong 0\;.
$$
Thus, the extension must be the trivial one and we conclude that 
$$
\widehat{H}^{\rm odd}(S^{2k+1};\hat{h})
\cong 
\ZZ/h\oplus \Omega^{\rm odd}_{d_H\mbox{-}{\rm cl}}(S^{2k+1})\;.
$$
\end{example}

\medskip
\paragraph{\bf Acknowledgement.}
The authors  thank the organizers and participants of the Geometric Analysis and  Topology Seminar at the
 Courant Institute at NYU for asking about twisting  Deligne cohomology, during a talk by H.S., which 
 encouraged the authors to finish this two-stage project, starting with \cite{GS5}.


\end{document}